\documentclass[3p,authoryear,11pt,fleqn,numbers,square,sort&compress]{elsarticle}
\usepackage{amssymb,amsbsy,amsmath,amsfonts,amssymb,amscd}
\usepackage[english]{babel}
\usepackage{color}
\usepackage{float}
\usepackage{amssymb,natbib}
\usepackage{graphicx}
\usepackage[latin1]{inputenc}
\usepackage{stmaryrd}
\usepackage{mathrsfs}
\usepackage{booktabs}
\usepackage{multirow,booktabs}
\usepackage{array}
\usepackage{subfig}
\usepackage{geometry}
\usepackage[ruled]{algorithm2e}
\usepackage{lscape}
\usepackage{titlesec}
\usepackage[hidelinks,colorlinks=true,citecolor=blue,linkcolor=red,urlcolor=green,pdfstartview=]{hyperref}
\usepackage{pifont}
\usepackage{lineno}
\usepackage{listings}
\usepackage{textcomp} 
\usepackage{appendix}
\usepackage[framed,numbered,autolinebreaks,useliterate]{mcode}    

\biboptions{comma,round}
\newtheorem{e-proposition}[theorem]{Proposition}

\newtheorem{e-definition}[theorem]{Definition\rm}

\def\og{\leavevmode\raise.3ex\hbox{$\scriptscriptstyle\langle\!\langle$~}}
\def\fg{\leavevmode\raise.3ex\hbox{~$\!\scriptscriptstyle\,\rangle\!\rangle$}}

\journal{arXiv}
\input{tcilatex}

\graphicspath{{fig/}}

\begin{document}

	\begin{frontmatter}
		
		
		
\title{Matrix-based implementation and GPU acceleration \\ of linearized ordinary state-based peridynamic models in MATLAB}
	\author[label1,label2]{Tao Ni\corref{cor1}}
	\ead{nitao_sklgp@cdut.edu.cn}		
\author[label2,label4]{Mirco Zaccariotto}	
	
	\author[label3]{Qizhi Zhu}

\author[label2,label4]{Ugo Galvanetto}

\cortext[cor1]{Corresponding author}
\address[label1]{State Key Laboratory of Geohazard Prevention and Geoenvironment Protection, Chengdu University of Technology, 610059 Chengdu, China}
\address[label2]{Industrial Engineering Department, University of Padova, via Venezia 1, Padova, 35131, Italy}
\address[label3]{College of Civil and Transportation Engineering, Hohai University, 210098, Nanjing, China}
\address[label4]{Center of Studies and Activities for Space (CISAS)-G. Colombo, via Venezia 15, Padova, 35131, Italy}
		\begin{abstract}
			Ordinary state-based peridynamic (OSB-PD) models have an unparalleled capability to simulate crack propagation phenomena in solids with arbitrary Poisson's ratio. However, their non-locality also leads to prohibitively high computational cost.  In this paper, a fast solution scheme for OSB-PD models based on matrix operation is introduced, with which, the graphics processing units (GPUs) are used to accelerate the computation. For the purpose of comparison and verification, a commonly used solution scheme based on loop operation is also presented. An in-house software is developed in MATLAB. Firstly, the vibration of a cantilever beam is solved for validating the loop- and matrix-based schemes by comparing the numerical solutions to those produced by a FEM software. Subsequently, two typical dynamic crack propagation problems are simulated to illustrate the effectiveness of the proposed schemes in solving dynamic fracture problems. Finally, the simulation of the Brokenshire torsion experiment is carried out by using the matrix-based scheme, and the similarity in the shapes of the experimental and numerical broken specimens further demonstrates the ability of the proposed approach to deal with 3D non-planar fracture problems. In addition, the speed-up of the matrix-based scheme with respect to the loop-based scheme and the performance of the GPU acceleration are investigated. The results emphasize the high computational efficiency of the matrix-based implementation scheme.
		\end{abstract}
		
		
	\end{frontmatter}
	
	\section{Introduction}
	
	
	Peridynamics (PD), firstly introduced by Silling in 2000 %
	\citep{silling2000reformulation}, is a non-local continuum theory
	based on integral-differential equations, which is an alternative new
	approach and has an unparalleled capability to simulate crack
	propagation in structures. Cracks can grow naturally without
	resorting to external crack growth criterion. In the past two decades,
	PD-based computational methods have made great progress in the simulation of crack propagation phenomena
	\citep{silling2005meshfree,lai2015peridynamics,wang2018three,cheng2019dynamic,zhang2019failure}%
	. The bond-based version of PD theory (BB-PD) was firstly presented in %
	\citep{silling2005meshfree}, and then was extended to its final version
	named state-based PD (SB-PD) in \citep{silling2007peridynamic}, which
	includes the ordinary and non-ordinary versions (OSB-PD and NOSB-PD).
	Several modified PD models were also proposed in %
	\citep{zhu2017peridynamic,wang20183,diana2019bond,zhang2019modified,liu2020new}. Although the
	PD-based numerical approaches have considerable advantages in solving crack
	propagation problems, they all share the same disadvantage of low
	computational efficiency due to the non-locality. 
	
	To overcome this shortcoming, there maybe two strategies, (i) reducing the computing costs or (ii) adopting parallel programming. In order to reduce the overall computing cost, several different methods \citep{galvanetto2016effective,han2016morphing,zaccariotto2017enhanced,zaccariotto2018coupling,fang2019method,Ni2019Coupling} have been proposed to couple PD-based models to models based on classical mechanics at continuous or discrete levels. In addition, the parallel programming of PD models has also been introduced in many papers. In \citep{lehoucq2008peridynamics}, the BB-PD models were implemented in the molecular dynamic package LAMMPS (Large-scale Atomic/Molecular Massively Parallel Simulator), where the calculation is paralleled via distributed-memory Message Passing Interface (MPI). $EMU$ is the first peridynamic code designed as a research code \citep{silling2005meshfree}, based on which, a scalable parallel code named $PDQ$ was developed and introduced in \citep{sakhavand2011parallel} for BB-PD models. An open-source computational peridynamics code named Peridigm, designed for massively-parallel multi-physics simulations and developed originally at Sandia National Laboratories, was firstly released in 2011 \citep{parks2012peridigm}, in which, the calculations can be performed in parallel with multi-cores. In \citep{lee2017parallel}, a parallel code for BB-PD model was developed in Fortran by using the Open Multi-Processing application programming interface (OpenMP), in which the PD model was coupled to the finite element method to reduce the computing costs. In \citep{li2020large}, the parallel implementation of BB-PD model on Sunway Taihulight Supercomputer was carried out. In \citep{diehl2020asynchronous}, an asynchronous and task-based implementation of PD models utilizing High Performance ParalleX (HPX) in C++ was introduced. In addition to the CPU-based (Central Processing Unit) parallelization schemes, GPU-based (Graphics Processing Unit) parallel programming is gaining more and more attention in the field of high-performance computing. In \citep{liu2012discretized}, a PD code was implemented using GPU for highly parallel computation via OpenACC, a programming standard designed to simplify parallel programming of heterogeneous CPU/GPU systems developed by Cray, CAPS, Nvidia and PGI. The Compute Unified Device Architecture (CUDA), which is a parallel architecture for NVIDIA GPUs, was used to implement PD model in \citep{diehl2012implementierung}. In \citep{zhang2016modeling}, OpenACC was adopted to accelerate the calculation of the internal PD force density on the CUDA-enabled GPU devices. Another alternative for the GPU-based implementation of PD models is the OpenCL, which is a framework for writing programs that run across heterogeneous platforms consisting of CPUs and GPUs, etc, a relevant work can be found in \citep{mossaiby2017opencl}. In many papers, the BB-PD model was preferred for parallel programming because of its lower computational cost. However, the classical BB-PD model has a limitation on the value of the Poisson's ratio, which has been removed in SB-PD model.
	Furthermore, the NOSB-PD model introduced the expression of stress into the
	equations but its numerical solutions are affected by zero energy mode \citep{silling2017stability}. Therefore, the OSB-PD model is the most
	convenient choice to describe the elastic deformation and fracture characteristics of solid materials with any Poisson's ratio \citep{bobaru2016handbook}. 
	
	It is a common opinion that vectorized calculation is usually faster than nested loop calculation, as confirmed in \citep{zhang2016modeling}. Inspired by that, we propose a matrix-based implementation scheme for linearized OSB-PD models and compare it with a non-parallel loop-based scheme. A modified explicit central difference time integration scheme proposed in \citep{taylor1989pronto} is adopted to obtain the dynamic solutions of the OSB-PD models, while the adaptive dynamic relaxation algorithm from \citep{Underwood1983dynamic} is adopted for the quasi-static solutions. An in-house software is developed in MATLAB and GPU acceleration is provided based on the built-in {\it Parallel Computing Toolbox}. Several typical numerical examples are carried out by using the developed software to investigate its capability and effectiveness in simulating elastodynamic and various crack propagation problems. Meanwhile, the performance of the loop- and matrix-based schemes, as well as the speed-up of the GPU acceleration, is also evaluated.
	
	
	This manuscript is organized as follows. The OSB-PD theory is briefly reviewed in $Sect. 2$. $Sect. 3$ describes the numerical implementation of the proposed solution schemes. In $Sect. 4$, several numerical examples are presented and discussed. Finally, $Sect. 5$ concludes the paper.
	
	\section{Ordinary state-based peridynamic theory}
	
	\subsection{Basic concepts}
	
	\begin{figure}[h]
		\centering
		\includegraphics[scale=0.4]{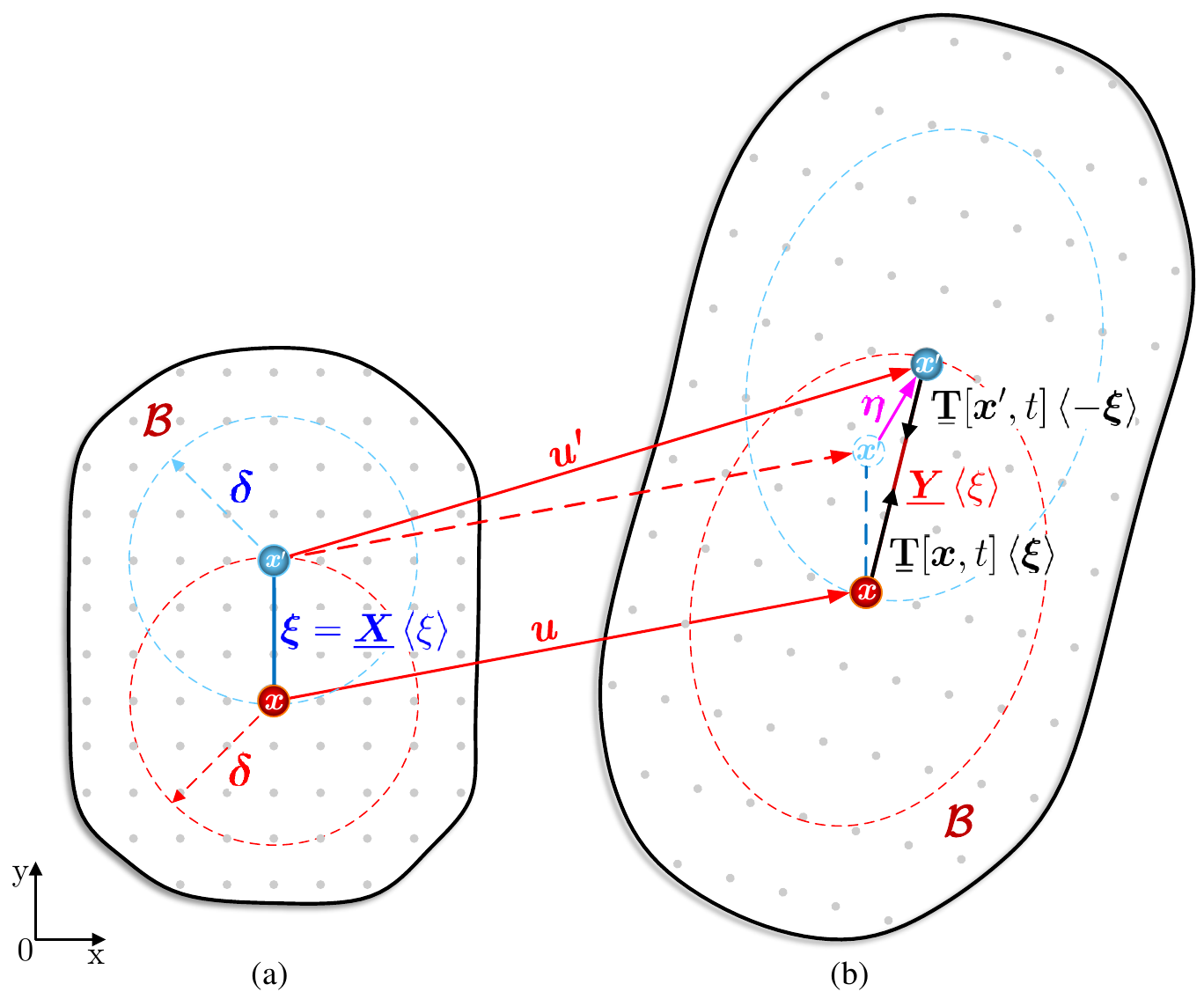}
		\caption{The positions of two points in the (a) initial and (b) deformed
			configurations.}
		\label{fig1}
	\end{figure}
	As shown in Fig. \ref{fig1}, a body $\left( \mathcal{B}\right) $ is modelled
	by OSB-PD, in which each material point interacts with all the other points
	around it within a region with a radius $\delta $ (called $horizon$) %
	\citep{silling2007peridynamic}. Assume that there are two points $\left( 
	\boldsymbol{x}\text{ and }\boldsymbol{x}^{\prime }\right) $ in the initial
	configuration of body $\mathcal{B}$, and the force density vector states of
	points $\boldsymbol{x}$ and $\boldsymbol{x}^{\prime }$ in the deformed
	configuration at time $t$\ are defined as \textbf{\b{T}}$\left[ \boldsymbol{x%
	},t\right] \left\langle \boldsymbol{\xi }\right\rangle $\ and \textbf{\b{T}}$%
	\left[ \boldsymbol{x}^{\prime },t\right] \left\langle \boldsymbol{-\xi }%
	\right\rangle $\ along the deformed bond, whose values can be different.
	Thus, the governing equation of motion of the ordinary state-based
	peridynamic system can be given as \citep{van2018objectivity}: 
	\begin{equation}
		\rho \boldsymbol{\ddot{u}}\left( \boldsymbol{x},t\right) =\int\nolimits_{%
			\mathcal{H}_{x}}\left\{ \text{\textbf{\b{T}}}\left[ \boldsymbol{x},t\right]
		\left\langle \boldsymbol{\xi }\right\rangle -\text{ \textbf{\b{T}}}\left[ 
		\boldsymbol{x}^{\prime },t\right] \left\langle \boldsymbol{-\xi }%
		\right\rangle \right\} dV_{x^{\prime }}+\boldsymbol{b}\left( \boldsymbol{x}%
		,t\right)  \label{1.1}
	\end{equation}%
	in which $\rho $ is the material mass density, $\boldsymbol{\ddot{u}}$ is
	the acceleration. $dV_{x^{\prime }}$ is the infinitesimal volume bundled
	with point $\boldsymbol{x}^{\prime }$. $\boldsymbol{b}$ is the force density
	applied by external force. $\mathcal{H}_{x}$ is the neighbourhood associated
	with the mass point $\boldsymbol{x}$, which is usually a circle in 2D and a
	sphere in 3D, mathematically defined as: $\mathcal{H}_{x}\left( \delta
	\right) =\left\{ x\in \mathcal{B}\text{:}\left\Vert \boldsymbol{\xi }%
	\right\Vert \leqslant \delta \right\} $.
	
	In ordinary state-based peridynamic theory, the concept of \textquotedblleft
	state\textquotedblright\ is introduced, the reference vector state and
	deformation vector state are defined as $\underline{\boldsymbol{X}}%
	\left\langle \boldsymbol{\xi }\right\rangle $ and $\underline{\boldsymbol{Y}}%
	\left\langle \boldsymbol{\xi }\right\rangle $, respectively, and expressed
	as:
	
	\begin{equation}
		\begin{array}{ccc}
			\underline{\boldsymbol{X}}\left\langle \boldsymbol{\xi }\right\rangle =%
			\boldsymbol{\xi } & , & \underline{\boldsymbol{Y}}\left\langle \boldsymbol{%
				\xi }\right\rangle =\boldsymbol{\xi +\eta }%
		\end{array}
		\label{1.2}
	\end{equation}
	where $\boldsymbol{\xi }$ is the initial relative position vector between
	points $\boldsymbol{x}$ and $\boldsymbol{x}^{\boldsymbol{\prime }}$ and $%
	\boldsymbol{\eta }$ is the relative displacement vector, which are defined
	as:%
	\begin{equation}
		\begin{array}{ccc}
			\boldsymbol{\xi =x}^{\boldsymbol{\prime }}-\boldsymbol{x} & , & \boldsymbol{%
				\eta }=\boldsymbol{u}^{\boldsymbol{\prime }}-\boldsymbol{u}%
		\end{array}
		\label{1.3}
	\end{equation}
	where $\boldsymbol{u}$ and $\boldsymbol{u}^{\boldsymbol{\prime }}\ $are the
	displacement vectors of points $\boldsymbol{x}$ and $\boldsymbol{x}^{%
		\boldsymbol{\prime }}$, respectively.
	
	The reference position scalar state and deformation scalar state are defined
	as:%
	\begin{equation}
		\begin{array}{ccc}
			\underline{x}=\left\Vert \underline{\boldsymbol{X}}\right\Vert & , & 
			\underline{y}=\left\Vert \underline{\boldsymbol{Y}}\right\Vert%
		\end{array}
		\label{1.5}
	\end{equation}
	where $\left\Vert \underline{\boldsymbol{X}}\right\Vert $ and $\left\Vert 
	\underline{\boldsymbol{Y}}\right\Vert $\ are the norms of $\underline{%
		\boldsymbol{X}}$ and $\underline{\boldsymbol{Y}}$, representing the lengths
	of the bond in its initial and deformed states, respectively.
	
	\subsection{Isotropic linear elastic solid material}
	
	In the deformed configuration of a linear elastic OSB-PD material (with
	Poisson' ratio $\nu $, bulk modulus $\kappa $\ and shear modulus $\mu $),
	the local deformation is described by two different components of volume
	dilatation value $\theta $ and the deviatoric extension state $\underline{e}%
	^{d}$,\ which are defined as %
	\citep{le2014two,sarego2016linearized,van2018objectivity}: 
	\begin{equation}
		\theta =\frac{A}{m}\int\nolimits_{\mathcal{H}_{x}} \underline{\mathit{w}}%
		\text{ }\underline{x}\text{ }\underline{e} dV_{x^{\prime }}  \label{1.6}
	\end{equation}%
	\begin{equation}
		\underline{e}^{d}=\underline{e}-\frac{\theta \underline{x}}{3}  \label{1.7}
	\end{equation}%
	in which, $A$ is a constant and given as: 
	\begin{equation}
		A=\left\{ 
		\begin{array}{ccc}
			3 & , & \text{3D} \\ 
			\frac{2\left( 1-2\nu\right) }{1-\nu } & , & \text{plane stress} \\ 
			2 & , & \text{plane strain}%
		\end{array}
		\right.  \label{1.8}
	\end{equation}
	%
	%
	$m$ is called the weighted volume and given as: 
	\begin{equation}
		m=\int\nolimits_{\mathcal{H}_{x}}\underline{\mathit{w}} \text{ }
		\underline{x}^{2}$d$V_{x^{\prime }}  \label{1.9}
	\end{equation}
	$\underline{\mathit{w}}$ is an influence function, whose forms have been
	summarized in \citep{Ni2019Coupling}. $\underline{e}$ is the extension
	scalar state for describing the longitudinal deformation of the bond,
	usually defined as:%
	\begin{equation}
		\underline{e}=\underline{y}-\underline{x}  \label{1.10}
	\end{equation}%
	The force density vector state $\mathbf{\b{T}}\left[ \boldsymbol{x},t\right]
	\left\langle \boldsymbol{\xi }\right\rangle $ is usually defined as: 
	\begin{equation}
		\text{\textbf{\b{T}}}\left[ \boldsymbol{x},t\right] \left\langle \boldsymbol{%
			\xi }\right\rangle =\underline{t}\cdot \underline{\boldsymbol{M}}%
		\left\langle \boldsymbol{\xi }\right\rangle  \label{1.11}
	\end{equation}%
	where $\underline{\boldsymbol{M}}\left\langle \boldsymbol{\xi }\right\rangle %
	$\ is a unit state in the direction of $\underline{\boldsymbol{Y}}$ defined
	as:%
	\begin{equation}
		\underline{\boldsymbol{M}}\left\langle \boldsymbol{\xi }\right\rangle =\frac{%
			\underline{\boldsymbol{Y}}}{\left\Vert \underline{\boldsymbol{Y}}\right\Vert 
		}  \label{1.13}
	\end{equation}%
	which has the approximation $\underline{\boldsymbol{M}}\left\langle 
	\boldsymbol{\xi }\right\rangle \approx {{\underline{\boldsymbol{X}}}}/{%
		\left\Vert \underline{\boldsymbol{X}}\right\Vert }$ under the hypothesis of
	infinitesimal deformations.
	
	$\underline{t}$ is called the force density scalar state, which can be
	obtained by \citep{van2018objectivity}: 
	\begin{equation}
		\begin{array}{ll}
			\underline{t}= & K\theta \frac{\underline{\mathit{w}}\text{ }\underline{x}}{m%
			}+G\underline{e}^{d}\frac{\underline{\mathit{w}}}{m}=\left( K-\frac{G}{3}%
			\right) \theta \frac{\underline{\mathit{w}}\text{ }\underline{x}}{m}+G%
			\underline{e}\frac{\underline{\mathit{w}}}{m}%
		\end{array}
		\label{1.14}
	\end{equation}%
	where $K$ and $G$\ are positive constants related to material parameters and
	given as:%
	\begin{equation}  \label{1.15}
		K=\left\{ 
		\begin{array}{ll}
			3\kappa & , \text{3D} \\ 
			2\kappa \frac{1-2\nu }{1-\nu }-\frac{2}{3}\mu \frac{\left( 1+\nu \right)
				\left( 1-3\nu \right) }{\left( 1-\nu \right) \left( 1-2\nu \right) } & ,%
			\text{ plane stress} \\ 
			2\kappa -\frac{2}{3}\mu & ,\text{ plane strain}%
		\end{array}%
		\right. \\
	\end{equation}
	
	\begin{equation}  \label{1.16}
		G=\left\{ 
		\begin{array}{ll}
			15\mu & , \text{3D} \\ 
			8\mu & ,\text{ plane stress and plane strain}%
		\end{array}%
		\right. \\
	\end{equation}
	
	\subsection{Failure criterion}
	
	In order to describe the failure and crack propagation in solids, failure
	criteria are essential in the PD models. The \textquotedblleft critical bond
	stretch\textquotedblright \ criterion first introduced in %
	\citep{silling2005meshfree} for BB-PD models is sometimes used for the
	OSB-PD models. However, different from that in BB-PD models, the
	deformation in OSB-PD models contains both the volumetric ($%
	\theta$) and deviatoric parts ($\underline{e}^d$). Therefore, the formulae
	of the criteria for OSB-PD models should be different from those for BB-PD
	models. Referring to \citep{zhang2018state}, a specific \textquotedblleft
	critical bond stretch\textquotedblright \ criterion is adopted for the
	OSB-PD models in this paper. The derivation of this criterion in 3D
	condition is explained in this section.
	
	\begin{figure}[h]
		\centering
		\includegraphics[scale=0.5]{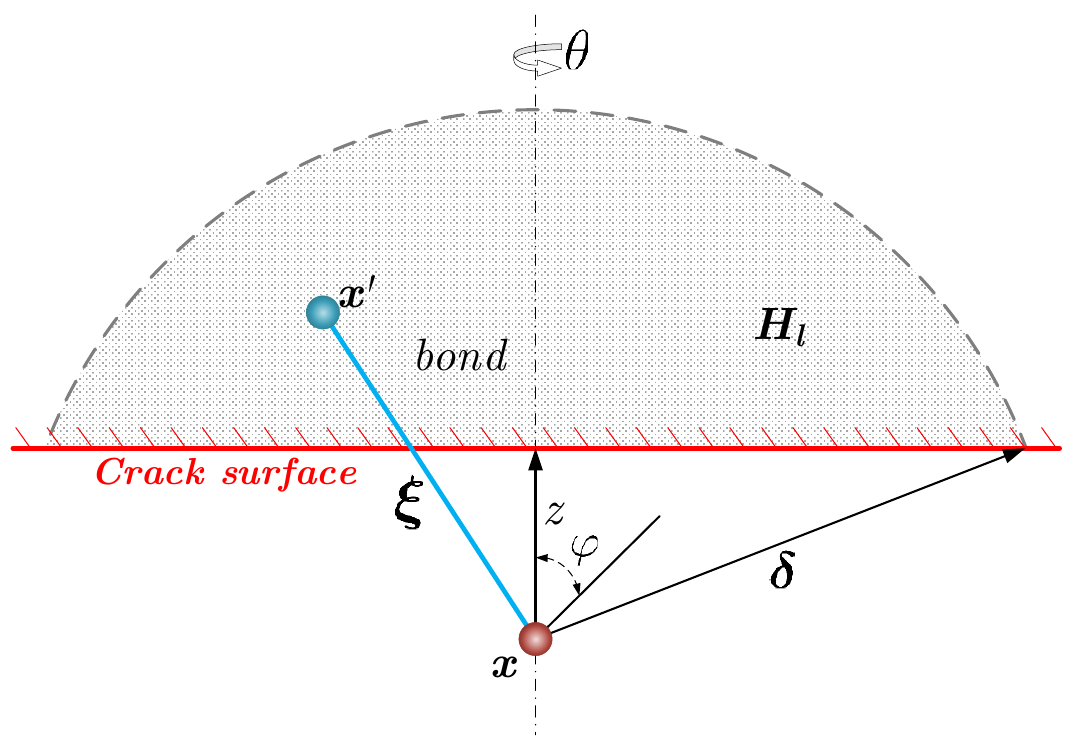}
		\caption{Schematic of a peridynamic domain crossed by a crack surface.}
		\label{fig2}
	\end{figure}
	
	Assuming that the critical value $s_c$ of all the bonds is the same. The
	bonds reaching the critical stretch value will be broken. Based on above
	notions, the stretch value of bond $\boldsymbol{\xi}$ can be defined as: 
	\begin{equation}
		s\langle\boldsymbol{\xi}\rangle=\frac{\underline{e}\langle\boldsymbol{\xi}%
			\rangle}{\underline{x}\langle\boldsymbol{\xi}\rangle}  \label{1.17}
	\end{equation}
	therefore, the extension scalar state can be rewritten as: 
	\begin{equation}
		\underline{e}\langle\boldsymbol{\xi}\rangle=s\langle\boldsymbol{\xi}\rangle 
		\underline{x}\langle\boldsymbol{\xi}\rangle  \label{1.18}
	\end{equation}
	
	As shown in Fig. \ref{fig2}, the neighbourhood of point $\boldsymbol{x}$ is
	crossed by a crack surface. $H_{l}$ represents the part removed by the crack
	from $\boldsymbol{x}$'s neighbourhood, and $\boldsymbol{x}^{\prime }$
	represents any point locating in $H_{l}$. The formation of cracked surface breaks the bond $\boldsymbol{\xi}$ and
	releases the strain energy stored in it. Then the work required to break all
	the bonds connecting point $\boldsymbol{x}$ to points in $H_l$ should be
	equal to the summation of the deformation energy stored in the broken bonds
	in their critical stretch condition. In an isotropic elastic OSB-PD
	material, the elastic strain energy density at point $\boldsymbol{x}$ in 3D
	condition can be expressed by a function of $\theta$ and $\underline{e}^{d}$%
	: 
	\begin{equation}
		W\left(\theta, \underline{e}^{d}\right)=\frac{1}{2}\kappa \theta^{2}+\frac{%
			15\mu}{2m} \int_{\mathcal{H}_{x}} \underline{w} \text{ } \underline{e}^{d} 
		\underline{e}^{d} \mathrm{~d} V_{x^{\prime}}  \label{1.19}
	\end{equation}
	
	As in \citep{zhang2018state,silling2010linearized}, substitution of Eqs.(\ref%
	{1.7}) and (\ref{1.9}) into Eq.(\ref{1.19}) rewrites the expression of elastic
	strain energy density at point $\boldsymbol{x}$ as: 
	\begin{equation}
		W(\theta, \underline{e})=\frac{1}{2}\left(\kappa-\frac{5}{3} \mu\right)
		\theta^{2}+\frac{15 \mu}{2m} \int_{\mathcal{H}_{x}} \underline{w} \text{ } 
		\underline{e} \text{ } \underline{e}\text{ } \mathrm{~d} V_{x^{\prime}}
		\label{1.20}
	\end{equation}
	
	According to Eq.(\ref{1.18}), the extension scalar state of the critically
	stretched bond can be expressed as $\underline{e}_c=s_c \underline{x}$.
	Thus, in the critical condition, the contribution of domain $H_l$ to the
	volume dilatation value at point $\boldsymbol{x}$ can be obtained by: 
	\begin{equation}
		\theta_{l}=\frac{3}{m} \int_{\mathcal{H}_{l}}\underline{w} \text{ } 
		\underline{x} \text{ } \underline{e}_{c} \mathrm{d} V_{x^{\prime}}=\frac{3}{m%
		} \int_{\mathcal{H}_{l}}\underline{w}\text{ } \underline{x}^{2} s_{c} 
		\mathrm{d} V_{x^{\prime}}=3 \varpi_{l} s_{c}  \label{1.21}
	\end{equation}
	where 
	\begin{equation}
		\varpi_{l}=\frac{m_{l}}{m} \quad \& \quad m_{l}=\int_{\mathcal{H}_{l}}%
		\underline{w} \underline{x}^{2} d V_{x^{\prime}}  \label{1.22}
	\end{equation}
	
	Based on the above formulae, the energy release by the broken bonds
	connecting points in domain $H_l$ to point $\boldsymbol{x}$ is: 
	\begin{equation}
		W\langle\boldsymbol{x}\rangle=\frac{1}{2}\left(\kappa-\frac{5}{3} \mu\right)
		\theta_{\boldsymbol{l}}^{2}+\frac{15 \mu}{2m} \int_{\mathcal{H}_{l}} 
		\underline{w} 
		\text{ }
		\underline{x}^{2} s_{c}^{2} \mathrm{~d} V_{x^{\prime}}=\frac{3}{2}\left[%
		\left(3\kappa-5\mu\right) \varpi_{l}^{2}+5 \mu \varpi_{l}\right] s_{c}^{2}
		\label{1.23}
	\end{equation}
	
	Recalling the above assumption, the critical energy $G_c$ released for per
	unit fracture area should conform to the following equation: 
	\begin{equation}
		G_{c}=2 \int_{0}^{\delta} W\langle\boldsymbol{x}\rangle \mathrm{d} z=2
		\int_{0}^{\delta}\frac{3}{2}\left[\left(3\kappa-5\mu\right) \varpi_{l}^{2}+5
		\mu \varpi_{l}\right] s_{c}^{2} \mathrm{~d} z=\left[\left(9\kappa-15
		\mu\right) \beta+15 \mu \beta^{\prime}\right] s_{c}^{2}  \label{1.24}
	\end{equation}
	in which, $\beta$ and $\beta^{\prime}$ are defined as: 
	\begin{equation}
		\beta=\int_{0}^{\delta} \varpi_{l}^{2} \mathrm{~d} z \quad \& \quad
		\beta^{\prime}=\int_{0}^{\delta} \varpi_{l} \mathrm{~d} z  \label{1.25}
	\end{equation}
	and their values depend on the discretization parameters and influence
	function.
	
	Based on Eq.(\ref{1.24}), the critical stretch value can be given
	accordingly as: 
	\begin{equation}
		s_{c}=\sqrt{\frac{G_{c}}{\left(9\kappa-15 \mu\right) \beta+15 \mu
				\beta^{\prime}}}  \label{1.26}
	\end{equation}
	
	In this paper, the influence function is taken as $w=1$, and the expression
	of $\beta$ and $\beta^{\prime}$ can be obtained as: 
	\begin{equation}
		\beta=\frac{125\delta}{1848} \quad, \quad \beta^{\prime}=\frac{5 \delta}{24
			\pi}  \label{1.27}
	\end{equation}
	
	Similarly, the formulae of the critical stretch value $s_{c}$ in 2D
	conditions can be obtained as: 
	\begin{equation}
		\begin{array}{ccc}
			s_{c}= \sqrt{\frac{G_{c}}{A^2\left(\kappa^{\prime}- \frac{8}{9} \mu\right)
					\beta+8 \mu \beta^{\prime}}} & , & \kappa^{\prime}=\left\{ 
			\begin{array}{ccc}
				\kappa+\frac{\mu}{9}(\frac{1+\nu}{1-2\nu})^2 & , & \text{plane stress} \\ 
				\kappa+\frac{\mu}{9} & , & \text{plane strain}%
			\end{array}%
			\right.%
		\end{array}\label{1.28}
	\end{equation}%
	in which, the definition of A is given in Eq.(\ref{1.8}), and the relevant
	expressions of $\beta $ and $\beta ^{\prime }$ are: 
	\begin{equation}
		\begin{array}{ccc}
			\beta = \frac{1087\delta }{1250\pi ^{2}} & , & \beta ^{\prime }=
			\frac{4\delta }{5\pi }%
		\end{array}\label{1.29}
	\end{equation}
	
	In order to indicate the connection status of the bonds, a scalar variable $%
	\mathit{\varrho }$ is defined as %
	\citep{zaccariotto2018coupling,ni2018peridynamic}: 
	\begin{equation}
		\underline{\mathit{\varrho }}\left\langle \boldsymbol{\xi }\right\rangle 
		\boldsymbol{=}\left\{ 
		\begin{array}{ccc}
			1 & , & \text{if }s\left\langle \boldsymbol{\xi }\right\rangle <s_{c} \\ 
			0 & , & \text{otherwise}%
		\end{array}%
		\right.  \label{1.30}
	\end{equation}%
	then the damage value $\varphi _{x}$ at point $\boldsymbol{x}$\ can be
	obtained by: 
	\begin{equation}
		\varphi _{x}=1-\frac{\int\nolimits_{\mathcal{H}_{x}}\underline{\mathit{w}}%
			\left\langle \boldsymbol{\xi }\right\rangle \text{ }\underline{\mathit{\
					\varrho }}\left\langle \boldsymbol{\xi }\right\rangle \text{d}V_{x^{\prime }}%
		}{\int\nolimits_{\mathcal{H}_{x}}\underline{\mathit{w}}\left\langle 
			\boldsymbol{\xi }\right\rangle \text{d}V_{x^{\prime }}}  \label{1.31}
	\end{equation}%
	in which $\varphi _{x}\in \left[ 0,1\right] $, and the cracks can be
	identified wherever $\varphi _{x}\geqslant 0.5$.
	
	\section{Discretization and numerical implementation}
	\subsection{Matrix-based discretization of the OSB-PD equations}
	After discretization, the peridynamic equation of motion of the current node 
	$\boldsymbol{x}_{i}$ at $n^{th}$ time increment is written in a form of
	summation:%
	\begin{equation}
		\rho \boldsymbol{\ddot{u}}_{i}^{n}=\sum\limits_{j=1}^{N_{H_{i}}}\left\{ 
		\text{\textbf{\b{T}}}\left[ \boldsymbol{x}_{i},n\right] \left\langle 
		\boldsymbol{\xi }_{ij}\right\rangle -\text{\textbf{\b{T}}}\left[ \boldsymbol{%
			x}_{j},n\right] \left\langle \boldsymbol{-\xi }_{ij}\right\rangle \right\}
		V_{j}+\boldsymbol{b}_{i}^{n}  \label{2.1}
	\end{equation}%
	where $N_{H_{i}}$ is is the number of family nodes of $\boldsymbol{x}_{i}$, $%
	\boldsymbol{x}_{j}$\ is $\boldsymbol{x}_{i}$'s family node, $V_{j}$\ is the
	volume of node $\boldsymbol{x}_{j}$. Eq.(\ref{2.1}) can also be rewritten as:
	\begin{equation}
		\rho V_{i}\boldsymbol{\ddot{u}}_{i}^{n}=\sum\limits_{j=1}^{N_{H_{i}}}\mathbf{%
			F}_{i}^{\boldsymbol{\xi }_{ij}}+\boldsymbol{b}_{i}^{n}V_{i}
	\end{equation}%
	in which, $V_{i}$ is the volume of node $\boldsymbol{x}_{i}$. $\mathbf{F}%
	_{i}^{\boldsymbol{\xi }_{ij}}$ is the internal force acting on
	node $\boldsymbol{x}_{i}$ through the deformed bond $\boldsymbol{\xi }_{ij}$%
	. Similarly, $\mathbf{F}
	_{j}^{\boldsymbol{\xi }_{ij}}$ is the force applied to
	node $\boldsymbol{x}_{j}$. Based on Eq.(\ref{1.14}), $\mathbf{F}_{i}^{\boldsymbol{\xi }_{ij}}$ and $\mathbf{F}_{j}^{\boldsymbol{\xi }_{ij}}$ can be
	expressed as:%
	\begin{equation}
		\begin{array}{c}
			\mathbf{F}_{i}^{\boldsymbol{\xi }_{ij}}=\left[ \left( K-\frac{G}{3}\right) 
			\tilde{\theta}_{ij}\underline{x}_{ij}+G\underline{\tilde{e}}_{ij}\right] 
			\underline{\mathit{w}}_{ij}\underline{\boldsymbol{M}}\left\langle 
			\boldsymbol{\xi }_{ij}\right\rangle V_{i}V_{j} \\
			\mathbf{F}_{j}^{\boldsymbol{\xi }_{ij}}=\left[ \left( K-\frac{G}{3}\right) 
			\tilde{\theta}_{ij}\underline{x}_{ij}+G\underline{\tilde{e}}_{ij}\right] 
			\underline{\mathit{w}}_{ij}\underline{\boldsymbol{M}}\left\langle 
			\boldsymbol{\xi }_{ji}\right\rangle V_{i}V_{j}
		\end{array}
		\label{2.2}
	\end{equation}%
	where $\tilde{\theta}_{ij}=\frac{\theta _{i}}{m_{i}}+\frac{\theta _{j}}{m_{j}%
	}$,\ $\underline{\tilde{e}}_{ij}=\frac{\underline{e}_{ij}}{m_{i}}+\frac{%
		\underline{e}_{ij}}{m_{j}}$\ and $\underline{\boldsymbol{M}}\left\langle 
	\boldsymbol{\xi }_{ij}\right\rangle =-\underline{\boldsymbol{M}}\left\langle 
	\boldsymbol{\xi }_{ji}\right\rangle $.
	
	According to Eq.(\ref{1.6}), the contributions of the deformed bond $%
	\boldsymbol{\xi }_{ij}$ to the volume dilatation values $\theta _{i}$ and $%
	\theta _{j}$ can be computed with:%
	\begin{equation}
		\theta _{i}^{\boldsymbol{\xi }_{ij}}=\frac{A\underline{\mathit{w}}_{ij}%
			\underline{x}_{ij}V_{j}}{m_{i}}\underline{e}_{ij}\text{, }\theta _{j}^{%
			\boldsymbol{\xi }_{ij}}=\frac{A\underline{\mathit{w}}_{ij}\underline{x}%
			_{ij}V_{i}}{m_{j}}\underline{e}_{ij}  \label{2.3}
	\end{equation}
	
	In this section, a 3D case is considered for the explanation of discretized equations. Supposing that $\boldsymbol{U}_i=[U_{i1},U_{i2},U_{i3}]$ and $\boldsymbol{U}_j=[U_{j1},U_{j2},U_{j3}]$ represent the displacement vectors of nodes $\boldsymbol{x}_{i}$ and $\boldsymbol{x}_{j}$, respectively. Then the value of $\underline{e}_{ij}$ can be computed by:
	\begin{equation}
		\underline{e}_{ij}=\left[\mathbf{C}^{\underline{e}}_{ij}\right]\left[
		\begin{array}{c}
			\boldsymbol{U}_i\\
			\boldsymbol{U}_j
		\end{array}\right] \label{2.4}
	\end{equation}
	
	In addition, Eqs. (\ref{2.2}) and (\ref{2.3}) can be rewritten respectively as following forms:
	\begin{equation}
		\left[ \begin{array}{c}
			\mathbf{F}_{i}^{\boldsymbol{\xi }_{ij}} \\ 
			\mathbf{F}_{j}^{\boldsymbol{\xi }_{ij}}%
		\end{array}%
		\right] =\left[\mathbf{K}_{ij}^{\theta }\right]\left[ 
		\begin{array}{c}
			\theta _{i} \\ 
			\theta _{j}%
		\end{array}%
		\right] +\left[\mathbf{K}_{ij}^{\underline{e}}\right]\underline{e}_{ij} \label{2.5}
	\end{equation}
	
	\begin{equation}
		\left[ \begin{array}{c}
			\theta _{i}^{\boldsymbol{\xi }_{ij}} \\
			\theta _{j}^{\boldsymbol{\xi }_{ij}}
		\end{array}\right]=\left[\mathbf{C}_{ij}^{\theta}\right]\underline{e}_{ij} \label{2.6}
	\end{equation}
	
	If the unit direction vector state of the bond $\boldsymbol{\xi }_{ij}$ is expressed as $\underline{\boldsymbol{M}}\left\langle 
	\boldsymbol{\xi }_{ij}\right\rangle=\left[M^{1}_{ij}, M^{2}_{ij}, M^{3}_{ij}\right]$,
	the matrices $\left[\mathbf{C}^{\underline{e}}_{ij}\right]$, $\left[\mathbf{K}_{ij}^{\theta }\right]$, $\left[\mathbf{K}_{ij}^{\underline{e}}\right]$ and $\left[\mathbf{C}_{ij}^{\underline{\theta}}\right]$ will be given as:
	\begin{equation}
		\left[\mathbf{C}^{\underline{e}}_{ij}\right]=\left[ 
		\begin{array}{cccccc}
			-M_{ij}^{1} & -M_{ij}^{2} & -M_{ij}^{3} & M_{ij}^{1} & M_{ij}^{2} & 
			M_{ij}^{3}%
		\end{array}%
		\right]\label{2.7}
	\end{equation}
	\begin{equation}
		\left[\mathbf{K}_{ij}^{\theta }\right]=\left( K-\frac{G}{3}\right) \underline{
			w}_{ij}\underline{x}_{ij}V_{i}
		V_{j}\left[ 
		\begin{array}{cccccc}
			\frac{M_{ij}^{1}}{\mathit{m}_{i}} & \frac{M_{ij}^{2}}{\mathit{%
					m}_{i}} & \frac{M_{ij}^{3}}{\mathit{m}_{i}} & -%
			\frac{M_{ij}^{1}}{\mathit{m}_{i}} & -\frac{M_{ij}^{2}}{\mathit{%
					m}_{i}} & -\frac{M_{ij}^{3}}{\mathit{m}_{i}} \\ 
			\frac{M_{ij}^{1}}{\mathit{m}_{j}} & \frac{M_{ij}^{2}}{\mathit{%
					m}_{j}} & \frac{M_{ij}^{3}}{\mathit{m}_{j}} & -%
			\frac{M_{ij}^{1}}{\mathit{m}_{j}} & -\frac{M_{ij}^{2}}{\mathit{%
					m}_{j}} & -\frac{M_{ij}^{3}}{\mathit{m}_{j}}%
		\end{array}%
		\right] ^{T} \label{2.8}
	\end{equation}
	\begin{equation}
		\left[\mathbf{K}_{ij}^{e}\right]=G\underline{w}_{ij}V_{i}V_{j}\left( \frac{1}{m_{i}}+\frac{%
			1}{m_{j}}\right) \left[ 
		\begin{array}{cccccc}
			M_{ij}^{1} & M_{ij}^{2} & M_{ij}^{3} & -M_{ij}^{1} & -M_{ij}^{2} & 
			-M_{ij}^{3}%
		\end{array}%
		\right] ^{T}  \label{2.9}
	\end{equation}
	\begin{equation}
	\left[\mathbf{C}_{ij}^{\theta }\right]=A\underline{w}_{ij}\underline{x}_{ij}\left[ 
		\begin{array}{cc}
			\frac{V_{j}}{m_{i}} & \frac{V_{i}}{m_{j}}%
		\end{array}%
		\right] ^{T}  \label{2.10}
	\end{equation}%
	
	Their forms in 2D conditions can be obtained by removing the terms related to the third coordinate component from the above formulae.
	
	\subsection{Time integration algorithm of the dynamic solution}
	The dynamic solution of the OSB-PD model is obtained by using a modified
		explicit central difference time integration scheme as discussed in %
		\citep{taylor1989pronto}, in which the velocities are integrated with a
		forward difference and the displacements with a backward difference. Then
		the velocity and displacement at the $\left( n+1\right) ^{th}$ time
		increment can be obtained as:%
		\begin{equation}
			\begin{array}{l}
				\boldsymbol{\dot{u}}^{n+1}=\boldsymbol{\dot{u}}^{n}+\Delta t\boldsymbol{%
					\ddot{u}}^{n} \\ 
				\boldsymbol{u}^{n+1}=\boldsymbol{u}^{n}+\Delta t\boldsymbol{\dot{u}}^{n+1}%
			\end{array}
			\label{2.11}
		\end{equation}%
		where $\boldsymbol{\ddot{u}}^{n}$ is the acceleration at $n^{th}$\ time
		increment and can be determined by using Newton's second law:%
		\begin{equation}
			\boldsymbol{\ddot{u}}^{n+1}=\mathbf{M}^{-1}\left( \mathbf{F}^{ext}-\mathbf{F}%
			^{int}\right)  \label{2.12}
		\end{equation}%
		in which $\mathbf{F}^{ext}$ and $\mathbf{F}^{int}$ are the external and
		internal force vectors, respectively, $\mathbf{M}$ is the diagonal mass
		matrix. $\Delta t$  is the constant time
		increment. An explicit method for the undamped system requires the use of a time step smaller than the
		critical time step for numerical stability. According to %
		\citep{zhou2016numerical}, the stable time increment for PD model can be
		defined as:%
		\begin{equation}
			\Delta t<\delta /c^{\prime }  \label{2.13}
		\end{equation}%
		where $c^{\prime }=\sqrt{\left( \lambda +2\mu \right) /\rho }$ is the
		dilatational wave speed and $\lambda $ and $\mu $ are the Lame's elastic
		constants of the material.
	\subsection{Implementation and GPU acceleration of the OSB-PD models in MATLAB} \label{NumericalIm}
\begin{figure}[h]
		\centering
		\includegraphics[scale=0.6]{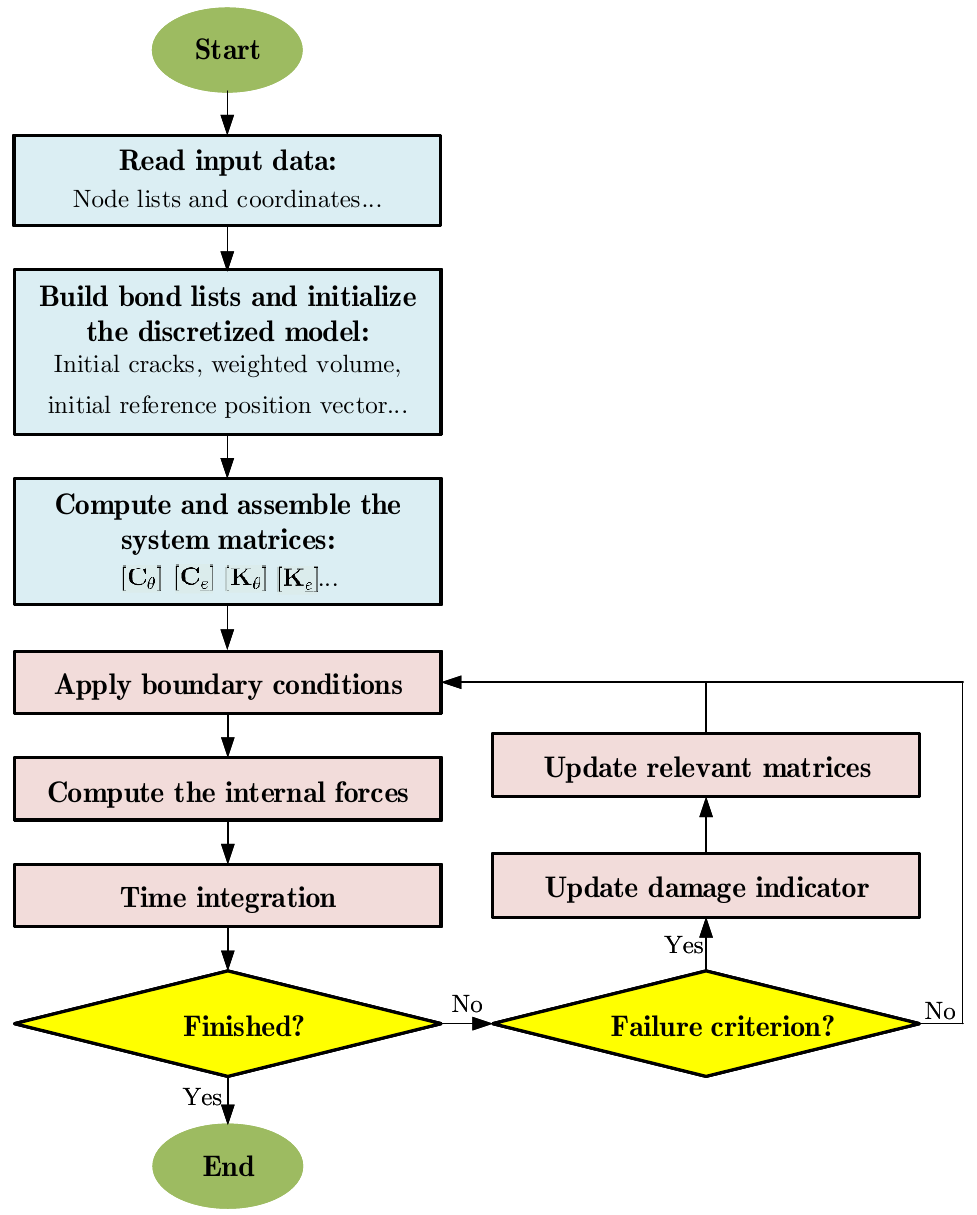}
		\caption{Roadmap for the matrix-based OSB-PD software implementation.}
		\label{fig3}
	\end{figure}

	MATLAB, an abbreviation of "matrix laboratory",  is a proprietary multi-paradigm programming language and computing environment, that has very high efficiency in matrix operations. For convenience, the matrix-based OSB-PD software is firstly developed in MATLAB. The roadmap for the implementation of the software is shown in Fig. \ref{fig3}. The rectangle blocks with light-grey background represent pre-processing sections, while the ones with pale orange background are solver sections. 
\begin{figure}[h]
	\centering
	\subfloat[]{\includegraphics[scale=0.8]{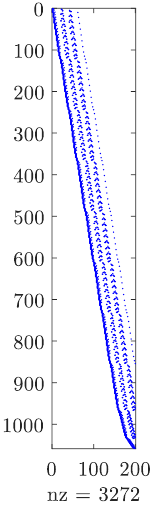}\label{fig3_1:sub1}}
	\hspace{0.5 in}
	\subfloat[]{\includegraphics[scale=0.8]{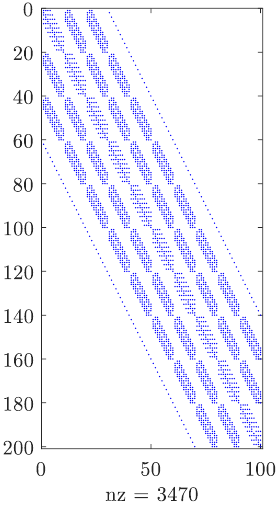}\label{fig3_1:sub2}}\\
	\subfloat[]{\includegraphics[scale=0.8]{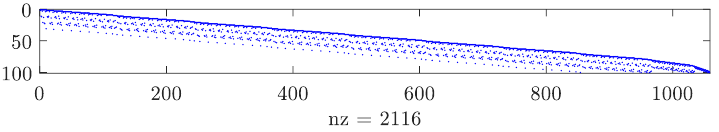}\label{fig3_1:sub3}}\\
	\subfloat[]{\includegraphics[scale=0.8]{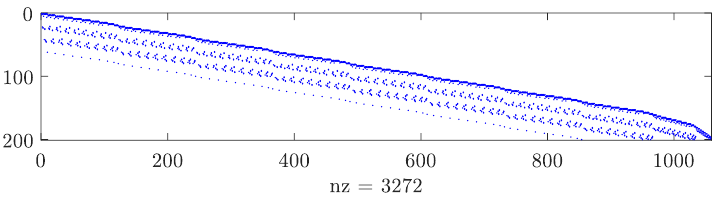}\label{fig3_1:sub4}}

	\caption{The shapes of the matrices (a) $[\mathbf{C}_{\underline{e}}]$ (b) $[\mathbf{K}_{\theta }] $ (c) $[\mathbf{C}_{\theta}]$ and (d) $[\mathbf{K}_{\underline{e}}]$ (\textquotedblleft nz\textquotedblright: number of non-zero entries).}
	\label{fig3_1}
\end{figure}

	$[\mathbf{K}_{\theta }] $, $[\mathbf{K}_{\underline{e}}]$, $[\mathbf{C}_{\underline{e}}] $ and $[ \mathbf{C}_{\theta}] $ are the global matrices of the system essential in the matrix-based OSB-PD software. Let us consider a 2D square discrete OSB-PD model with $10\times10$ nodes, when the $m$-ratio is taken equal $3$, the number of the PD bonds is $1058$. The shape of the system matrices with minimum bandwidths are shown in Figs. \ref{fig3_1:sub1} to \ref{fig3_1:sub4}. Obviously, these matrices are sparse and can be obtained by assembling the matrices of each bond. In the software, the system matrices are built and stored by using MATLAB function \textquotedblleft $sparse$\textquotedblright.
	
	 Considering a 3D discrete peridynamic model consisting of $N_{n}$ nodes and $N_{b}$ bonds, then the initialization of the variables required for the solver in MATLAB environment can be done as shown in appendix \ref{APX1}, and an example of the implementation of the solver sections is shown in appendix \ref{APX2}. For the purpose of comparison, the implementation scheme based on loop operation is also developed and an example of MATLAB function for the calculation of PD forces is shown in appendix \ref{APX3}. In order to make the loop operation in MATLAB have a computational efficiency close to that of the C language, the introduced function is pre-compiled by C-mex.

In a simulation adopting explicit iterative solution algorithms, the pre-processing sections will be executed only once, while the solver sections will be executed multiple times. Therefore, the key to improve the overall calculation speed is to improve the execution efficiency of the solver sections. In MATLAB, the {\it Parallel Computing Toolbox} is provided to perform parallel computations on multicore computers, GPUs, and computer clusters, which can significantly improve the efficiency of program execution when solving computationally and data-intensive problems. In the solver sections of the matrix-based implementation scheme, most of the calculations are vectorized. Therefore, the acceleration of the simulation can be easily achieved by copying the vectors and matrices used in the solver to GPUs, and this operation is performed by using MATLAB function \textquotedblleft $gpuArray$\textquotedblright.

\begin{figure}[h]
	\centering
	\subfloat[]{\includegraphics[scale=0.8]{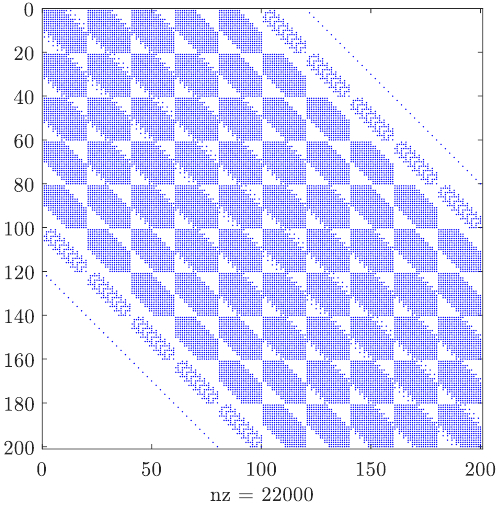}\label{fig3_2:sub1}}\hspace{0.5 in}
	\subfloat[]{\includegraphics[scale=0.8]{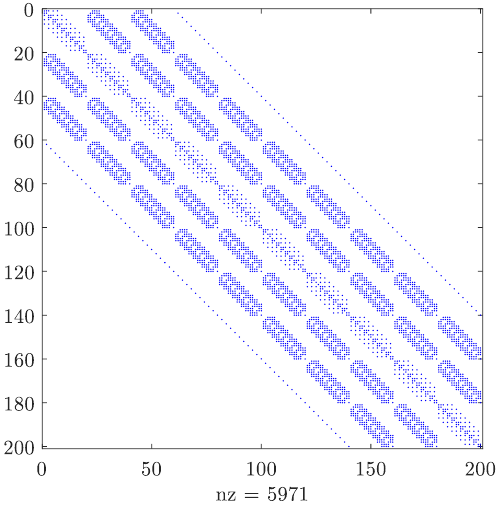}\label{fig3_2:sub2}}\\
	\subfloat[]{\includegraphics[scale=0.8]{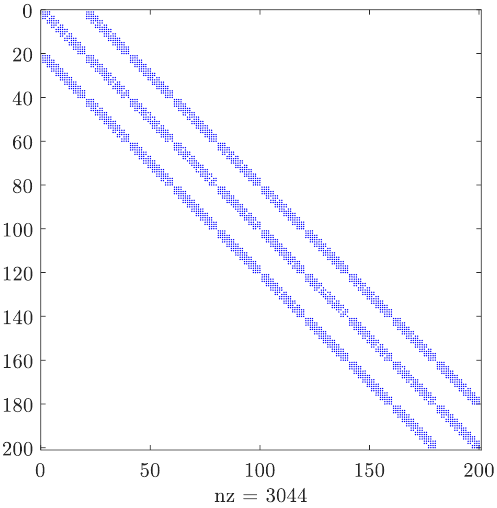}\label{fig3_2:sub3}}
	\caption{The shapes of the global stiffness matrices generated by (a) OSB-PD model (b) BB-PD model (c) FEM model (\textquotedblleft nz\textquotedblright: number of non-zero entries). Note that the three models are characterised by the same number of degrees of freedom.}
	\label{fig3_2}
\end{figure}

\textbf{Remark 1}. According to Eqs. (\ref{2.4}), (\ref{2.5}) and (\ref{2.6}), the following expression can be obtained:
\begin{equation}
	\left[\mathbf{F}^{int}\right] =\left[\mathbf{K}^{v}\right]\left[\mathbf{U}\right] +\left[\mathbf{K}^{d}\right]\left[\mathbf{U}\right]=\left(\left[\mathbf{K}^{v}\right]+\left[\mathbf{K}^{d}\right]\right)\left[\mathbf{U}\right]  \label{2.14}
\end{equation}
where $\left[\mathbf{K}^{v}\right]=\left[\mathbf{K}_{\theta}\right]\cdot \left[\mathbf{C}_{\theta }\right]\cdot\left[\mathbf{C}_{\underline{e}}\right]$ and $\left[\mathbf{K}_{d}\right]=\left[\mathbf{K}_{e}\right]\cdot\left[\mathbf{C}_{\underline{e}}\right]$. It is obvious that $\left(\left[\mathbf{K}^{v}\right]+\left[\mathbf{K}^{d}\right]\right)$ is the global stiffness matrix of the discretized equilibrium equations in the displacement-force form. The shape of the stiffness matrix corresponding to the matrices described in Figs. \ref{fig3_1:sub1} to \ref{fig3_1:sub4} is explained in Fig. \ref{fig3_2:sub1}, while the shape of the stiffness matrices obtained by BB-PD model and FEM model are shown in Figs. \ref{fig3_2:sub2} and \ref{fig3_2:sub3}. These figures indicate that the  bandwidth of the OSB-PD stiffness matrix is greater than that of BB-PD model, and is significantly greater than that of FEM model. From this perspective, it will be both memory- and time-intensive to solve the OSB-PD model by using an implicit algorithm. Thus, coupling PD-based models to FE models will be convenient to improve computational efficiency \citep{galvanetto2016effective,zaccariotto2018coupling,Ni2019Coupling,ni2019static,sun2019superposition,ni2020hybrid,dong2020stability,liu2021coupling}, moreover an iterative solver is suggested in the simulation.

\clearpage

\section{Numerical examples}
	
In this section, several examples are presented to demonstrate the effectiveness and the performance of the proposed schemes. All the cases are discretized using uniform grids and the $m$-ratio is always taken as $3$.
	
The in-house software for the simulations is developed in MATLAB 2019b, and executed on a Desktop computer with Intel$^{{\tiny ®}}$ Xeon$^{{\tiny ®}}$ CPU E5-1650 V3, 3.50 GHz processor and 64 GB of RAM. The GPU acceleration is performed on a NVIDIA GeForce GTX 1080 Ti. The peak Double-Precision Performance of the CPU and GPU are around 336 GFLOPs and 11067.4 GFLOPs, respectively. All cases tested in this section are solved in double-precision, thus, the theoretical peak acceleration ratio is around 32.9.

\subsection{Example 1: vibration of a cantilever beam}	
\begin{figure}[h]
	\centering
	\includegraphics[scale=0.7]{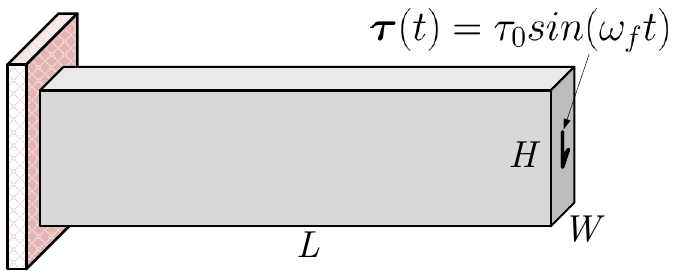}
	\caption{Geometry of the cantilever beam and boundary conditions used in the simulation.}
	\label{fig6}
\end{figure}

In this section, an elastodynamic problem of the vibration of a cantilever beam described in \citep{mossaiby2017opencl} is investigated without considering any crack propagations. As shown in Fig. \ref{fig6}, the geometric parameters of the beam specimen are taken as: $L=48m$, $H=12m$ and $W=5m$. The left end of the beam is constrained, while the right end is set free and driven by a downward uniformly distributed shear stress $\boldsymbol{\tau}(t)=\tau_{0}sin(\omega_{f} t) $, where the magnitude and the frequency of the load are taken as $\tau_{0}=13.89 Pa$ and $\omega_{f}=27 rad/s$, respectively. The mechanical parameters of the beam specimen are given as: Young modulus: $E=30 MPa$, Poisson's ratio: $\nu=0.25$, and mass density $\rho =1 kg/m^{3}$. 

\begin{table}[h]
	\caption{Discretization parameters of the cantilever beam.}
	\label{tab1.1}\centering {\scriptsize 
		\begin{tabular}{p{3cm}<{\centering}p{1cm}<{\centering}p{1cm}<{\centering}p{1cm}<{\centering}p{1cm}<{\centering}p{1cm}<{\centering}p{1cm}<{\centering}}
			\toprule Case & 1 & 2 & 3 & 4 & 5 & 6 \\ \hline
			Horizon radius ($\delta$) $[m]$ & 2.4 & 1.2 & 0.6 & 0.3 & 0.15& 0.075 \\ \hline
			Number of nodes & 1024 & 3844 & 14884 & 58564 & 232324& 925444\\ 
			\bottomrule &  & &  &  &  &
		\end{tabular}%
	}
\end{table}

\begin{figure}[h]
	\centering
    \includegraphics[scale=0.75]{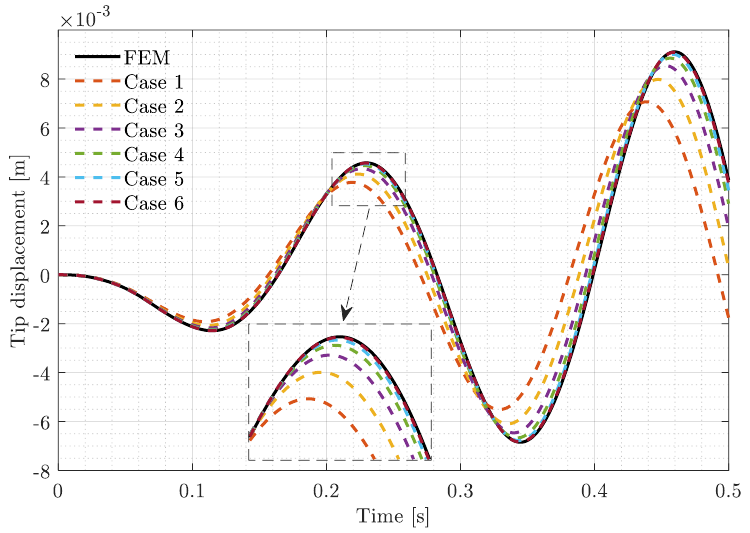}
	\caption{Variation of the vertical displacement at beam tip versus time in 2D simulations.}
	\label{fig7}
\end{figure}

The problem is firstly solved in plane stress condition. As listed in Tab. \ref{tab1.1}, six cases with different horizon sizes are carried out. The total simulation duration is $0.5s$ and a fixed time step of $\Delta t=20 \mu s$ is adopted. All the cases in Tab. \ref{tab1.1} are solved by using the matrix-based scheme, the variations of the tip displacement versus time during the simulations are recorded and plotted in Fig. \ref{fig7}. As shown in the magnifying box of Fig. \ref{fig7}, as the grid size decreases, the difference between the OSB-PD solution and FEM solution gradually decreases. 


\begin{figure}[h]
	\centering
	\subfloat[Computing times of different cases.]{\includegraphics[scale=0.75]{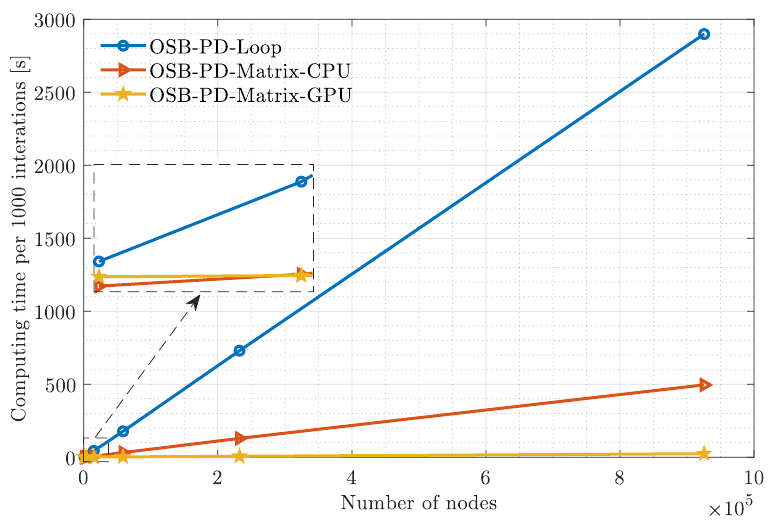}\label{fig8:sub1}}\\
	\subfloat[Speed-up in different cases.]{\includegraphics[scale=0.75]{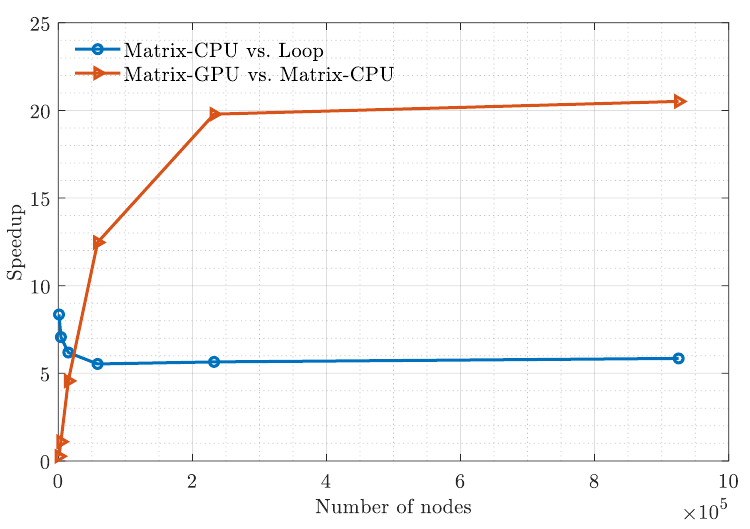}\label{fig8:sub2}}
	\caption{Comparison of the performance of different implementation schemes.}
	\label{fig8}
\end{figure}

In order to investigate the efficiency of different schemes, the cases in Tab. \ref{tab1.1} are carried out by using both the loop- and matrix-based schemes, and the matrix-based scheme is performed both on CPU and GPU. The computing times spent for $1000$ iterations by using different schemes are plotted in Fig. \ref{fig8:sub1}. Fig. \ref{fig8:sub2} shows the speed-up of the matrix-based scheme to loop-based scheme and the GPU acceleration ratio to the matrix-based scheme. It is worth noting that in the cases with a smaller amount of calculation, the acceleration of the GPU with respect to the matrix-based scheme is poor. As the amount of calculation increases, the acceleration effect becomes gradually more significant, and the maximum speed-up ratio is greater than $20$. Moreover, in the cases with a smaller amount of calculation, the speed-up ratio of the matrix-based scheme with respect to the loop-based scheme is greater that $6$. As the amount of calculation increases, the speed-up ratio is around $6$, which is usually affected by different factors, such as the CPU's turbo frequency parameters and the memory size of the computer etc... 

\begin{figure}[h]
	\centering
	\includegraphics[scale=0.75]{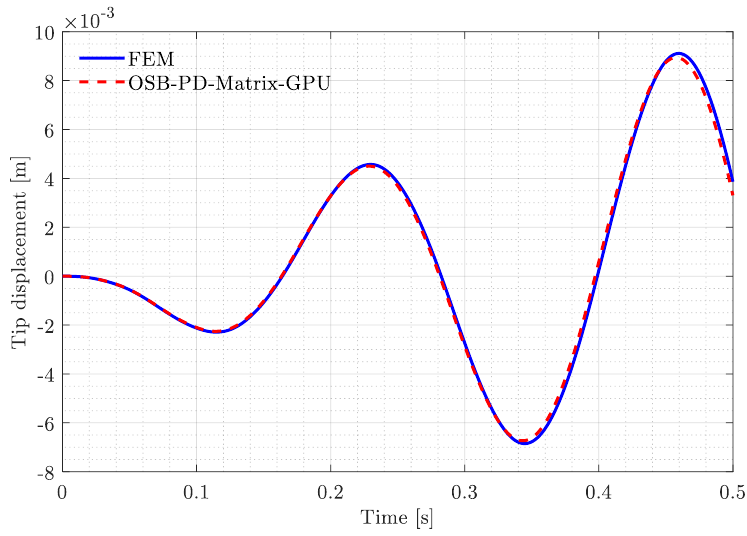}
	\caption{Variation of the vertical displacement at beam tip versus time in 3D simulation.}
	\label{fig9}
\end{figure}

In addition, the described problem is also solved in 3D by using the matrix-based scheme. A horizon radius of $\delta=0.6m$ is adopted for the discretization, and the grid size is $\Delta x=\delta/3=0.2m$, which results in $386984$ nodes. The total simulation duration is $0.5s$ and a fixed time step of $\Delta t=20 \mu s$ is adopted for the time integration. The computing time for $1000$ iterations running on CPU and GPU is $1005.7s$ and $39.413s$, respectively. The corresponding GPU acceleration ratio is around $25.517$. The variation of tip displacement versus time is compared with FEM solution and plotted in Fig. \ref{fig9}, showing that the solution obtained by the matrix-based scheme is in good agreement with that of FEM.
\newpage
\subsection{Example 2: pre-cracked plate subjected to traction}
\begin{figure}[h]
	\centering
	\includegraphics[scale=0.475]{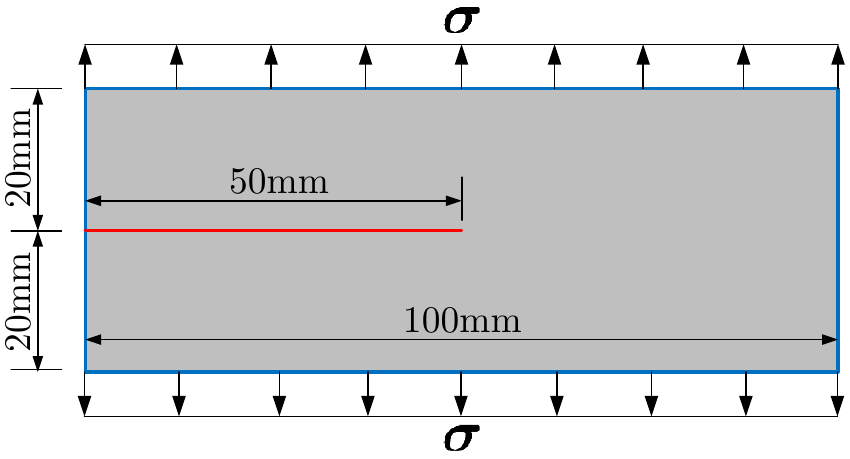}
	\caption{Geometry and boundary conditions of the pre-cracked plate subjected
		to traction.}
	\label{fig4_1}
\end{figure}

In this section, a pre-cracked plate subjected to traction is simulated to test the performance of the matrix-based scheme in simulating dynamic crack branching phenomena.  The geometry and boundary conditions are shown
in Fig. \ref{fig4_1}. The traction is kept constant during the whole
duration of the simulation, and two cases with different traction loads are considered, case 1: $\sigma =
20 MPa$ and case 2: $\sigma = 40 MPa$. The material parameters used in %
\citep{shojaei2018adaptive} are adopted, Young modulus: $E=72 GPa$, mass
density $\rho =2235 kg/m^{3}$, Poisson's ratio: $\nu=1/3$ (plane stress
condition), and fracture energy density: $G_{0}=380 J/m^2$. 

A horizon radius of $\delta=0.75mm$ is used for the discretization, and the corresponding grid size is adopted as $\Delta x=\delta/3=0.25mm$. In the discrete model, the total number of
the nodes is 64962. The simulation time durations of cases 1 and 2 are
chosen as $40 \mu s $ and $30 \mu s$, respectively, and a fixed time step of 
$\Delta t=25 ns$ is used for the time integration. %

\begin{figure}[h!]
	\centering  
	\subfloat[Case 1 obtained by the proposed scheme]{\includegraphics[width=0.445\textwidth]{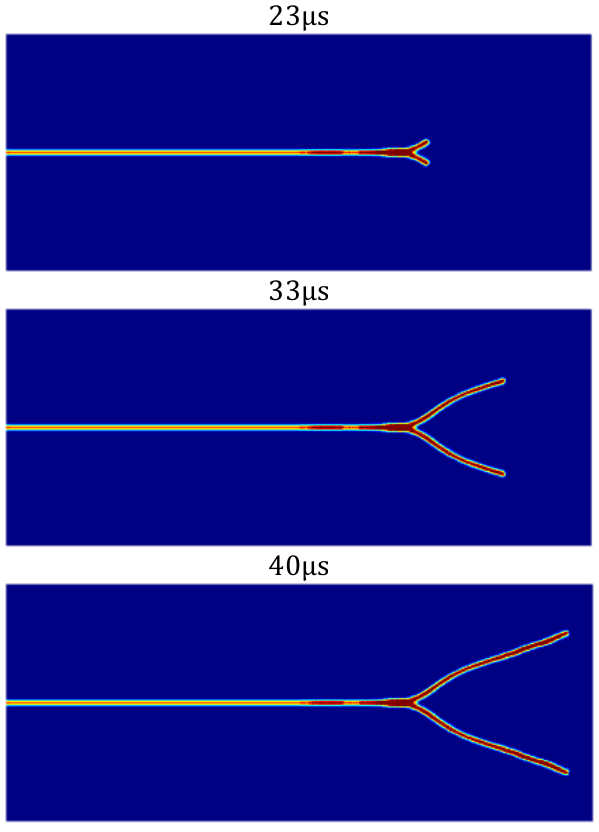}\label{fig4_42:sub1}}
	\subfloat[Case 1 with BB-PD \citep{shojaei2018adaptive}]{\includegraphics[width=0.45\textwidth]{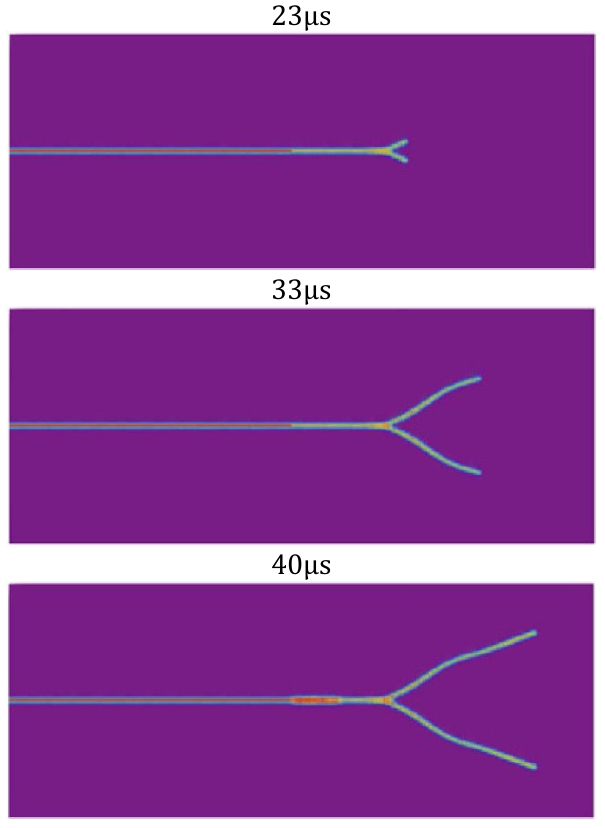}\label{fig4_42:sub2}}
	\caption{Damage contours at $23\mu s$, $33\mu s$ and $40\mu s$ of the case 1 in example 2.}
	\label{fig4_42}
\end{figure}

\begin{figure}[h!]
	\centering  
	\subfloat[Case 2 obtained by the proposed scheme]{\includegraphics[width=0.45\textwidth]{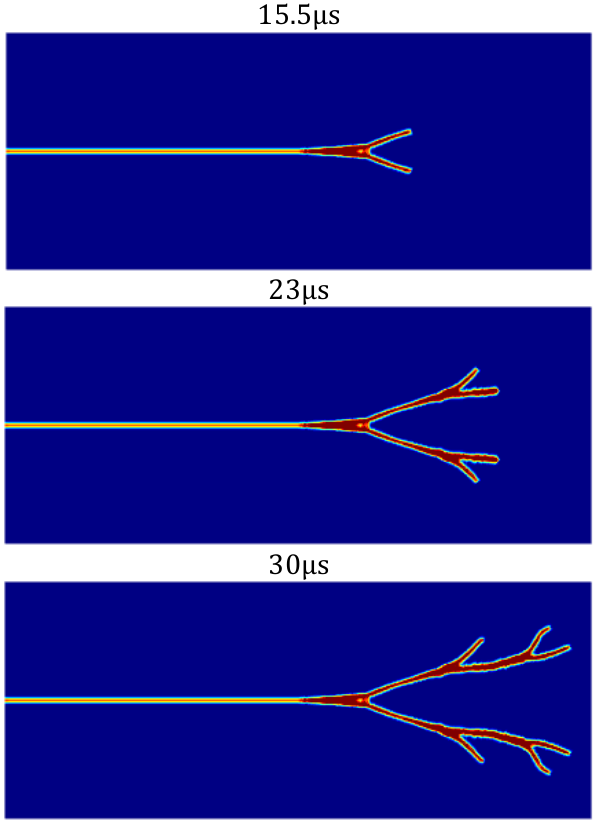}\label{fig4_43:sub1}}
	\subfloat[Case 2 with BB-PD \citep{shojaei2018adaptive}]{\includegraphics[width=0.44\textwidth]{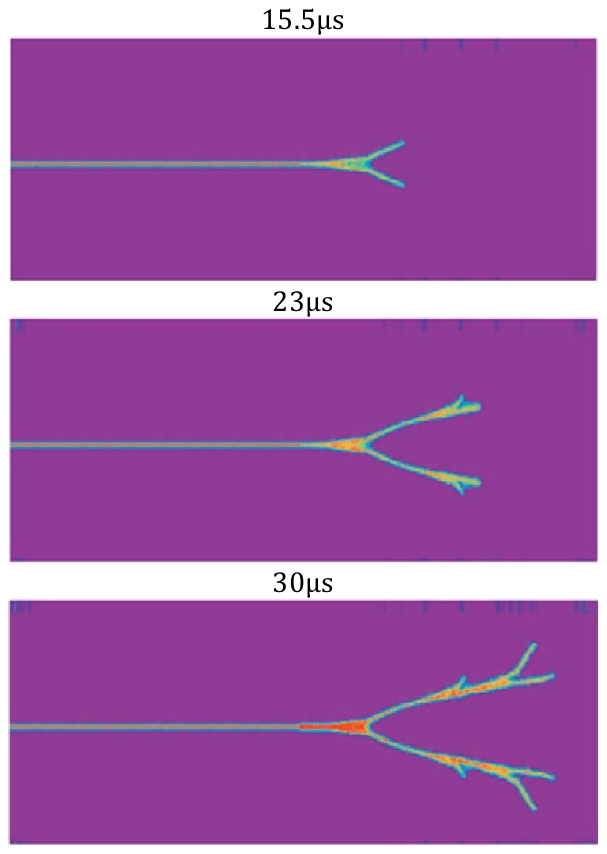}\label{fig4_43:sub2}}
	\caption{Damage contours at at $15.5\mu s$, $23\mu s$ and $30\mu s$ of the case 2 in example 2.}
	\label{fig4_43}
\end{figure}

The OSB-PD solutions at several selected time instants of cases 1 and 2 obtained by the matrix-based scheme are shown in Figs. \ref{fig4_42:sub1} and \ref{fig4_43:sub1}, while Figs. %
\ref{fig4_42:sub2} and \ref{fig4_43:sub2} are the corresponding numerical
results from \citep{shojaei2018adaptive}. The numerical results obtained by the proposed scheme and those from \citep{shojaei2018adaptive} share similar characteristics in the location and number of bifurcations. In \citep{shojaei2018adaptive}, the bond-based peridynamic model is adopted, and with the same critical fracture energy release rate, the critical bond stretch values computed by Eqs.(\ref{1.26}) and (\ref{1.28}) are smaller than those in \citep{shojaei2018adaptive}. Therefore, it is readily comprehensible that the simulated crack patterns are different from those in \citep{shojaei2018adaptive}.


In addition, the models are also solved by the loop-based scheme for comparing the computational efficiencies, and the computing times required by the loop-based and matrix-based schemes are shown in Tab. \ref{tab2}. In cases 1 and 2, the speed-up ratios of the matrix-based solution procedure executed on GPU with respect to that executed on CPU are around $5.13$ and $4.57$, respectively. The speed-up ratios of the matrix-based scheme with respect to the loop-based scheme, both executed on CPU, are around $7.16$ and $6.9$, respectively. 

\begin{table}[h!]
	\caption{Computing costs of example 2.}
	\label{tab2}\centering {\scriptsize 
		\begin{tabular}{p{4cm}<{\centering}p{3cm}<{\centering}p{3cm}<{\centering}p{3cm}<{\centering}}
			\toprule Methods & Loop-based scheme & Matrix-based scheme (CPU) & Matrix-based scheme (GPU) \\ \hline
			Computing time of case 1 [s] & 350.81 & 48.967 & 9.54\\ \hline
			Computing time of case 2 [s] & 264.4 & 38.29 & 8.37\\ 
			\bottomrule &  & &
		\end{tabular}%
	}
\end{table}

\newpage
\subsection{Example 3: Kalthoff-Winkler's experiment}
	
\begin{figure}[h]
		
		\centering
		\includegraphics[width=4cm,height=6cm]{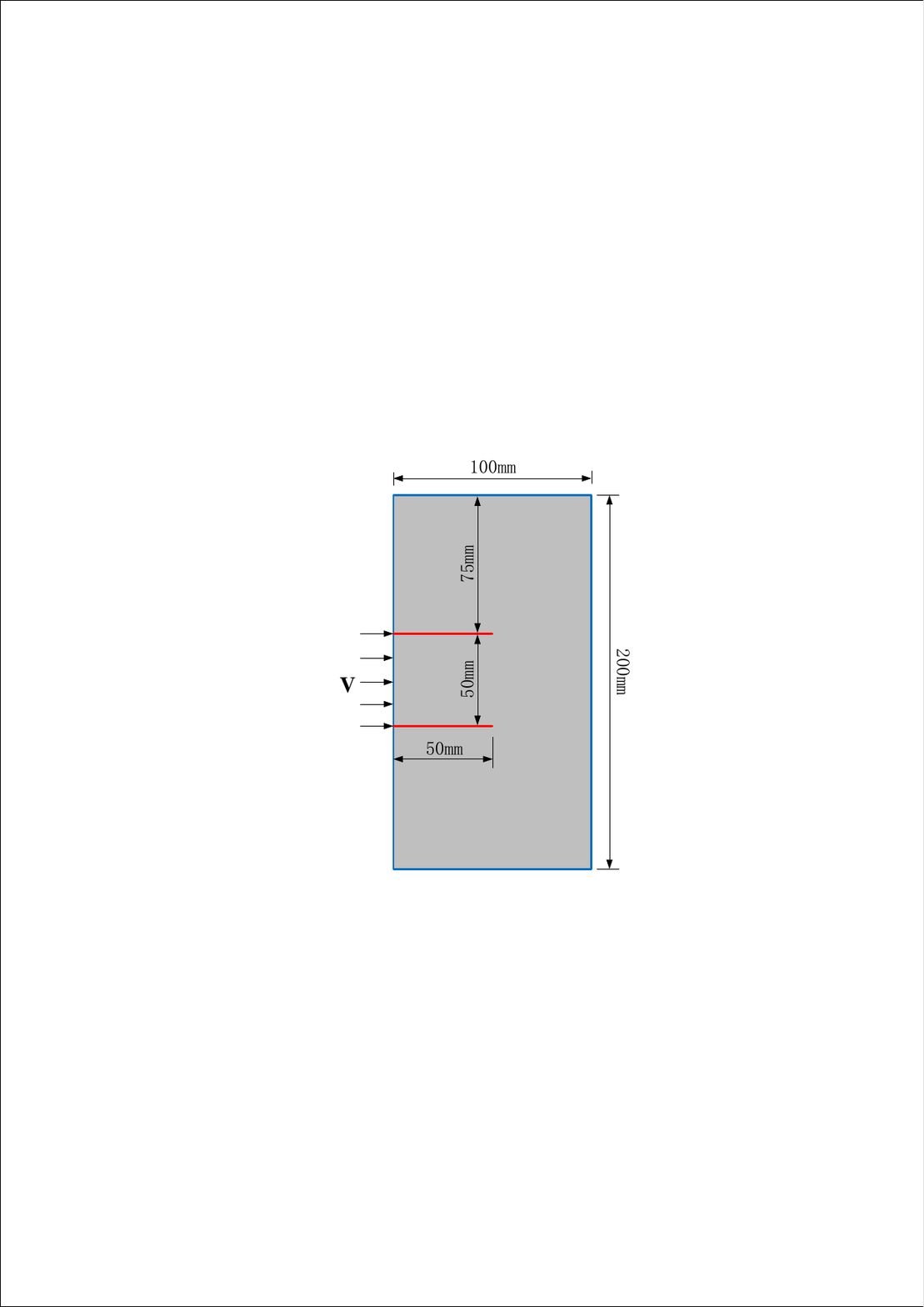}
		\caption{Geometry and boundary conditions of the Kalthoff-Winkler's
			experiment.}
		\label{fig4_5}
	\end{figure}

	The third example is the Kalthoff-Winkler's experiment reported in \citep{kalthoff2000modes}, which is a typical dynamic fracture setting
	used for the validation of numerical methods. The geometry and boundary conditions
	are explained in Fig. \ref{fig4_5}. The material parameters are given as,
	Young modulus: $E=190 GPa$, mass density $\rho =7800 kg/m^{3}$, Poisson's
	ratio: $\nu=0.25$ (plane strain condition), and fracture energy density: $%
	G_{0}=6.9\times10^{4} J/m^2$.
	
	 The horizon radius is adopted as $\delta =1.5mm$,
	and the corresponding grid size is $\Delta x=\delta/3 =0.5mm$. The
	total number of the nodes used for discretization is 80802. An initial
	horizontal velocity of $V=22 m/s$ is applied to the first three layers of
	nodes between the notches shown in Fig. \ref{fig4_5} and remains constant
	during the simulation. The total simulation duration is $t=100 \mu s$ and a
	fixed time step of $\Delta t=20 ns $ is chosen for the time integration. 
	
\begin{figure}[h!]
		\centering
		\subfloat[Result obtained by the proposed scheme]{%
			\includegraphics[width=11.5cm,height=6.75cm]{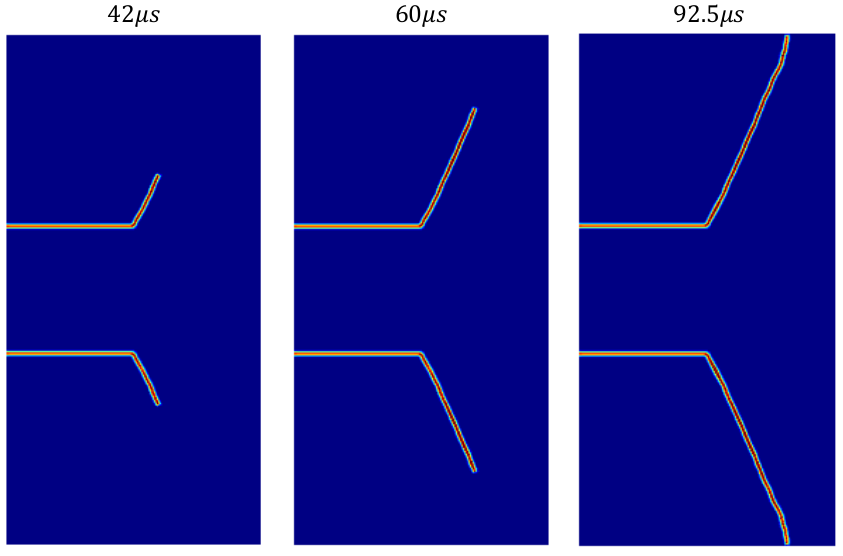}\label{fig4_44:sub1}} \\
		\subfloat[Result from \citep{ren2016dual}]{\includegraphics[width=11.5cm,height=6.75cm]{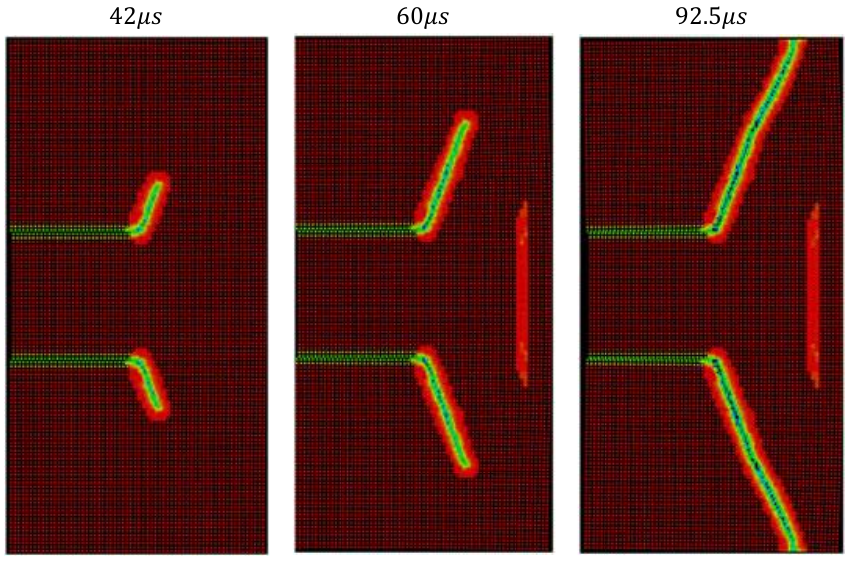}\label{fig4_44:sub2}}
		\caption{Damage contours at $42\mu s$, $60\mu s$ and $92.5\mu s$ in example 3.}
		\label{fig4_44}
	\end{figure}

	The crack pattern obtained by the OSB-PD model solved with the matrix-based scheme is shown in Fig. \ref%
	{fig4_44:sub1}, while Fig. \ref{fig4_44:sub2} is the numerical result from \citep{ren2016dual}. The comparison of the crack patterns
	shown in Figs. \ref{fig4_44:sub1} and \ref{fig4_44:sub2} indicates the accuracy of the presented
	matrix-based solution scheme. In order to compare the computational efficiencies, the model is also solved by the loop-based scheme. The computing times spent by the loop-based and matrix-based schemes are shown in the Tab. \ref{tab1}. The calculation speed of the matrix-based
	scheme executed on CPU is higher than that of the loop-based scheme on CPU, and the speed-up ratio is around $5.6$, while the ratio of the GPU acceleration is around  $6.18$.
	\begin{table}[h!]
		\caption{Computing costs of example 3.}
		\label{tab1}\centering {\scriptsize 
		\begin{tabular}{p{4cm}<{\centering}p{3cm}<{\centering}p{3cm}<{\centering}p{3cm}<{\centering}}
	\toprule Methods & Loop-based scheme & Matrix-based scheme (CPU) & Matrix-based scheme (GPU) \\ \hline
	Computing time [s] & 1392.62 & 248.9 & 40.29\\ 
	\bottomrule &  & &
\end{tabular}}	
\end{table}

\subsection{Example 4: Brokenshire torsion experiment}
\begin{figure}[h]
	\centering
	\includegraphics[scale=0.4]{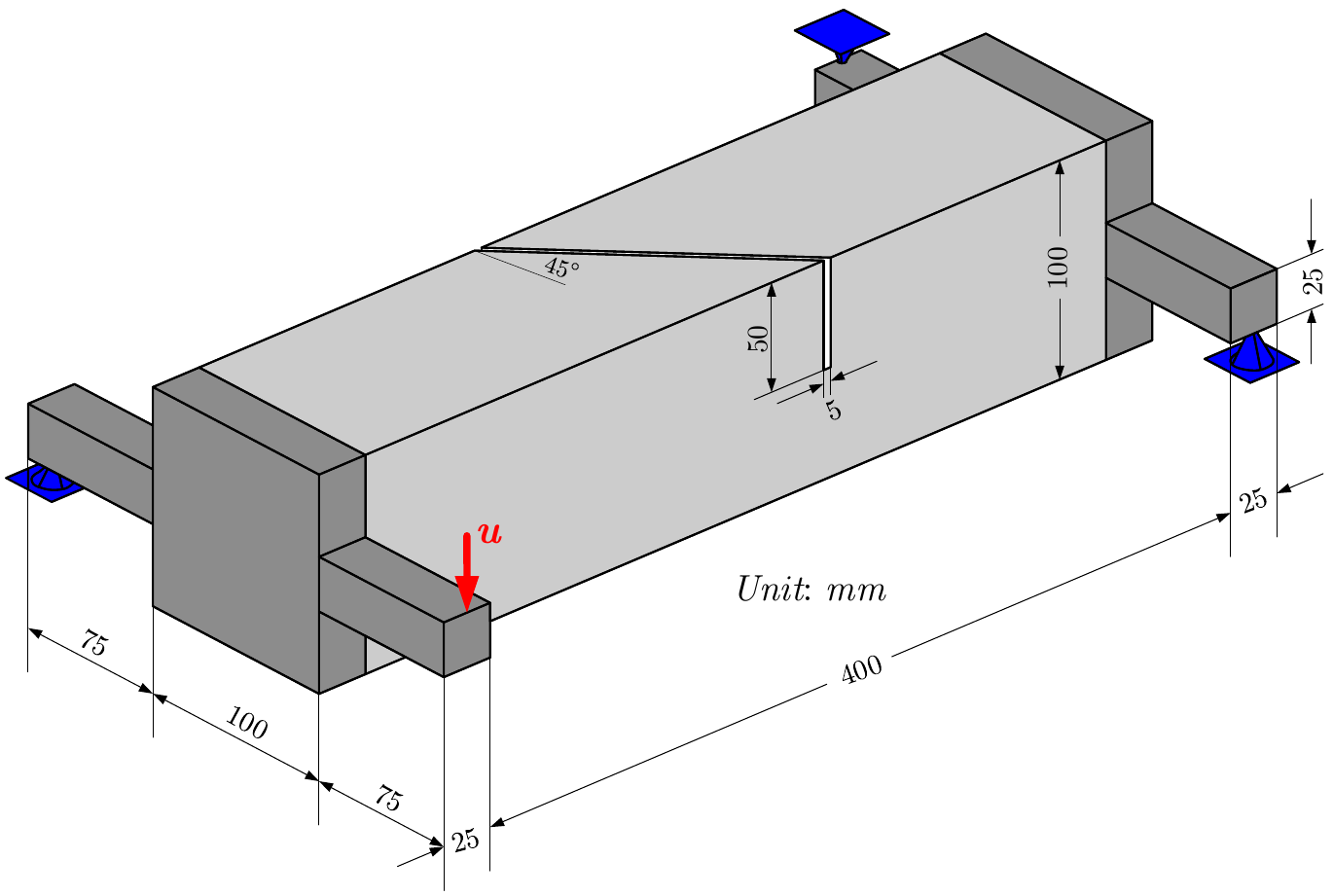}
	\caption{Geometry and boundary conditions of the Brokenshire torsion experiment.}
	\label{fig5_1}
\end{figure}

The purpose of the last example, the Brokenshire torsion experiment reported in \citep{brokenshire1995study}, is to illustrate the capability of the OSB-PD model implemented by using the proposed matrix-based scheme to simulate non-planar 3D crack propagation problems. The main geometrical parameters of the prismatic specimen and the boundary conditions of the test are shown in Fig.\ref{fig5_1}.
The material parameters of the specimen are taken as, Young
modulus: $E=35GPa$, Poisson's ratio: $\nu =0.2$, and $G_{0}=80J/m^{2}$ \citep{jefferson2004three}.  
A horizon radius of $\delta=7.5mm$ is adopted for discretization, and the corresponding grid size is $\Delta x=2.5mm$. The number of the nodes in the discretized model is 314961. 

Differently from the time integration algorithm adopted in the previous cases, the model is solved by using an adaptive dynamic relaxation algorithm presented in \citep{Underwood1983dynamic,kilic2010adaptive,Ni2019Coupling} for its quasi-static solution. A gradually increasing vertical downward displacement is applied to the tip of the loading arm (see Fig.\ref{fig5_1}) with a fix increment of $\Delta u=5\times10^{-8}m$, and $50000$ iterations are performed. 

\begin{figure}[h]
	\centering
	\subfloat[]{%
		\includegraphics[scale=0.475]{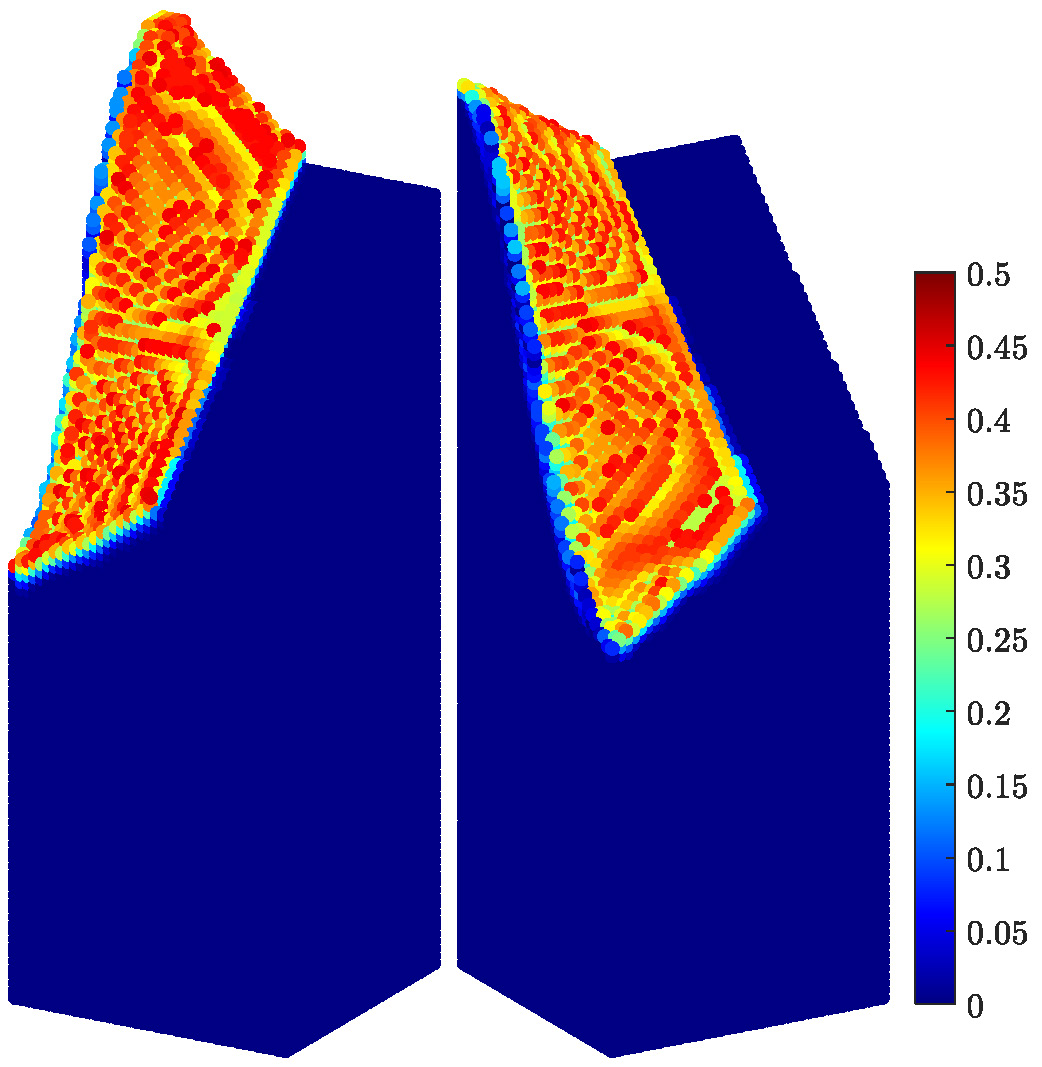}\label{fig5_2:sub1}} \hspace{0.25in}
	\subfloat[]{\includegraphics[scale=0.6]{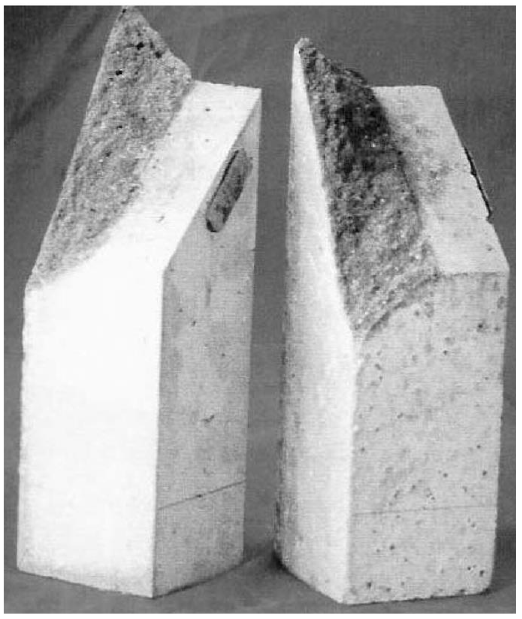}%
		\label{fig5_2:sub2}} 

	\caption{The (a) numerically and (b) experimentally \citep{jefferson2004three} observed broken specimens.}
	\label{fig5_2}
\end{figure}

The simulated broken specimen is shown in \ref{fig5_2:sub1}, while the experimentally observed broken specimen reported in \citep{jefferson2004three} is shown in Fig. \ref{fig5_2:sub2}. The numerically observed fracture surfaces are very similar to those observed experimentally. The computing times of the simulation carried out on the CPU and GPU are $50664.294s$ and $2489.646s$, respectively, where a speed-up ratio of about $20.35$ is achieved.
\newpage
\section{Conclusions}
	
	A matrix-based implementation for the OSB-PD models is proposed in this paper. An in-house OSB-PD software is developed in MATLAB and GPU acceleration is easily achieved. In order to validate the matrix-based scheme, a commonly used scheme in the meshless implementation of PD models based on loop operation is also introduced, and the loop calculation is compiled into the mex functions to achieve in MATLAB a computational efficiency comparable with C language.
	
	Firstly, the vibration of a cantilever beam is solved by using the matrix-based scheme both in plane stress and 3D conditions. Several cases with different horizon radius values are adopted for the 2D simulations, the numerical results show that the difference between OSB-PD and FEM solutions decreases with horizon value reduction. Subsequently, two different dynamic crack propagation problems are solved by using the proposed scheme. The numerical results produced by the matrix-based scheme are generally very similar to those in \citep{shojaei2018adaptive} and \citep{ren2016dual}, which indicates the accuracy of the solution scheme presented in this paper. Finally, the Brokenshire torsion experiment, which is a typical 3D non-planar fracture example, is simulated by using the matrix-based scheme, and the comparison of the fracture features observed in the experimental and numerical broken specimens further validates the effectiveness of the proposed scheme. 
	
	In addition, the speed-up of the matrix-based scheme with respect to loop-based scheme executed on CPU and the speed-up ratio in the GPU acceleration of the matrix-based scheme are also investigated. The results show that when the overall calculation amount of the model is small, the speed-up of the matrix-based scheme with respect to the loop-based scheme is more significant, whereas the performance of GPU acceleration is poor. With the increase of the calculation amount, the speed-up ratio of the matrix-based scheme to the loop-based scheme gradually decreases until it stabilizes at around $6$, while the performance of GPU acceleration gradually improves, and the speed-up ratio can reach around $20$. Note that, with different computing devices, the investigated speed-up ratios will be somewhat different, but the computational efficiency of the matrix-based scheme is certainly higher than that of loop-based scheme.
	
In general, the implementation scheme based on matrix operation in MATLAB can greatly improve the efficiency and simplify the GPU parallel programming of the OSB-PD models in the simulation of crack propagation problems, which minimizes programming effort, maximizes performance and paves the way for the application of the OSB-PD models to solve engineering problems with a large number of degrees of freedom. It should be noted that the proposed matrix-based scheme can also be implemented in other languages with efficient execution libraries for matrix operations (such as Python and Julia) maintaining a high efficiency.
	
\clearpage	
\begin{appendices}
\section{Initialization of the variables required for the solver executed on CPU} \label{APX1}
\begin{lstlisting}
%Array of the displacement components of each node.	
U=zeros(3*Nn,1);
%Array of the velocity components of each node.	
VEL=zeros(3*Nn,1);
%Array of the velocity components of each node in previous iteration.	
VELold=zeros(3*Nn,1);
%Array of the acceleration components of each node.	
ACC=zeros(3*Nn,1);
%Array of the peridynamic force components of each node.	
PDF=zeros(3*Nn,1);
%Array of the external force components of each node.
EXF=zeros(3*Nn,1);
%Array of the residual force components of each node.
RF=zeros(3*Nn,1);
%Array of the dilatation values of each node.
Theta=zeros(Nn,1);
%Array of the extension scalar state value of each bond.
E=zeros(Nb,1);
\end{lstlisting}
\newpage
\section{An example of implementation of the solver based on matrix operations}\label{APX2}
\begin{lstlisting}
for it=1:IntTimeStep
	%===========Block for applying Boundary Conditions============% 
	...
	%====Block for the computation of Peridynamic force Array====% 
	%Calculate the extension scalar values of each bond.
	E=CE*U;
	%Compute the dilatation values of each node.
	Theta=CTH*E;
	%Calculate the peridynamic internal force components of each node.
	PDF=KTH*Theta+KE*E;
	%=======Update the displacement components of each node=======%
	%Calculate the residual force components of each node.
	RF=EXF+PDF;
	%Calculate the acceleration components of each node.
	ACC=RF./M;
	%Calculate the current velocity components of each node.
	VEL=VELold+ACC*Dt;
	%Calculate the current displacement components of each node.
	U(FreeDofs)=U(FreeDofs)+Dt*VEL(FreeDofs); 
	%Store the current velocity components of each node.
	VELold=VEL;
	%============Block for applying failure ctriterion============% 
	...  
	%=================Block for the result output=================% 
	...       
end
\end{lstlisting}
\newpage
\section{An example of peridynamic forces calculation function based on loop operation}\label{APX3}
\begin{lstlisting}
function [PDF]=PDForceComputLoop(U,NodeVol,BondList,m,xi,Xi,CoIndex)
	%NodeVol: Array of the volumes of each node.
	%BondList: Array of the connected index of each bond.
	%m: Array of the weight volumes of each node.
	%xi: Array of the initial reference position scalar state values of each bond.
	%Xi: Array of the initial reference position Array state components of each bond.
	%CoIndex: Tag array of bond connection status: 1 represents connection, 0 represents broken.
	%Array of the dilatation values of each node.
	Theta=zeros(Nn,1);
	%Array of the peridynamic force components of each node.
	PDF=zeros(Nn*2,1);    
	%Array of the extension scalar state value of each bond.
	E=zeros(Nb,1);
	%Array of the initial reference position scalar state values of each bond.
	Eta=zeros(3*Nb,1);
	%Array of the initial reference position Array state components of each bond.
	eta=zeros(Nb,1);  
	%Calculate the extension scalar values of each bond. 
	Eta(1:3:end)=Xi(1:3:end)+U(3*BondList(:,2)-2)-U(3*BondList(:,1)-2);
	Eta(2:3:end)=Xi(2:3:end)+U(3*BondList(:,2)-1)-U(3*BondList(:,1)-1);
	Eta(3:3:end)=Xi(3:3:end)+U(3*BondList(:,2))-U(3*BondList(:,1));
	eta=sqrt(Eta(1:3:end).^2+Eta(2:3:end).^2+Eta(3:3:end).^2);
	E=eta-xi;
	E(CoIndex==0)=0;
	%Calculate the dilatation values of each node. 
	for i=1:Nb
	    index1=BondList(i,1);index2=BondList(i,2);
	    Theta(index1)=Theta(index1)+A*xi(i)*E(i)*CoIndex(i)*NodeVol(index2);
	    Theta(index2)=Theta(index1)+A*xi(i)*E(i)*CoIndex(i)*NodeVol(index1);
	end 
	Theta=Theta./m;
	%Calculate the peridynamic force components of each node.
	for i=1:Nb
	    index1=BondList(i,1);index2=BondList(i,2);
	    Mx=Xi(3*i-2)/xi(i);My=Xi(3*i-1)/xi(i);Mz=Xi(3*i)/xi(i);
	    DT1=Theta(index1)/m(index1);DT2=Theta(index2)/m(index2);
	    DE1=E(i)/m(index1);DE2=E(i)/m(index1);
	    FT=(K-G/3)*(DT1+DT2)*xi(i)*CoIndex(i)*NodeVol(index1)*NodeVol(index2);
	    FE=G*(DE1+DE2)*CoIndex(i)*NodeVol(index1)*NodeVol(index2);  
	    PDF(3*index1-2)=PDF(3*index1-2)+(FT+FE)*Mx;
	    PDF(3*index1-1)=PDF(3*index1-1)+(FT+FE)*My;
	    PDF(3*index1)=PDF(3*index1)+(FT+FE)*Mz;
	    PDF(3*index2-2)=PDF(3*index2-2)-(FT+FE)*Mx;
	    PDF(3*index2-1)=PDF(3*index2-1)-(FT+FE)*My;
	    PDF(3*index2)=PDF(3*index2)-(FT+FE)*Mz;        
	end
end	
\end{lstlisting}	
\end{appendices}

%
%
\clearpage
\section*{Acknowledgements}

This research is financially supported by the National Natural Science Foundation of China (Grant No. 42207226); Natural Science Foundation of Sichuan Province (Grant No. 2023NSFSC0808); State Key Laboratory of Geohazard Prevention and Geoenvironment Protection Independent Research Project SKLGP2021Z026.

The authors would like to acknowledge the support they received from MIUR
under the research project PRIN2017-DEVISU and from University of Padua
under the research project BIRD2020 NR.202824/20.

\end{document}